\documentclass[10pt,journal,twocolumn]{IEEEtran} 
\usepackage{amsmath,amsfonts}
\usepackage{amssymb}
\usepackage{graphicx}
\usepackage[colorlinks=true, allcolors=blue]{hyperref}
\usepackage{bbm}
\usepackage{caption,subcaption}
\usepackage{verbatim}
\usepackage{cite,enumerate}
\usepackage{srcltx,cases}
\usepackage[mathscr]{eucal}
\usepackage{multirow}
\usepackage{comment}

\newtheorem{theorem}{Theorem}
\newtheorem{lemma}{Lemma}
\newtheorem{proposition}{Proposition}

\newtheorem{remark}{\textit{Remark}}

\newtheorem{example}{Example}

\title{The Skyline Process:\\
Quantifying Sky Visibility\\ in 3D Urban Environments}
\author{Junse Lee, Fran\c{c}ois Baccelli\thanks{Junse Lee is with School of AI Convergence, Sungshin Women's University, Seoul, South Korea (e-mail: junselee@sungshin.ac.kr).}
\thanks{Fran\c{c}ois Baccelli is with INRIA/ENS and Telecom Paris (e-mail francois.baccelli@ens.fr).}}

\begin{document}
\maketitle

\begin{abstract}
Non-terrestrial networks (NTNs) are considered a promising technology for seamless, universal communication in the 6G era. However, signals from NTN elements to ground users are often blocked by high-rise buildings in dense urban environments. To quantify this blocking effect, in this paper, we propose a novel analytical framework by modeling the location of buildings as a 3D skyline process based on stochastic geometry and we derive closed-form expressions for the distribution of the blockage elevation angles. Furthermore, we extend our analysis to include spatially correlated blockage effects and analyze the spectral properties of the Skyline process using power spectral density (PSD) and autocorrelation function (ACF). Based on these theoretical findings, we present numerical results that provide insights into the design of LEO satellite networks by computing the mean number of visible satellites and the outage probability. We provide a decorrelation angle that serves as a useful threshold for obtaining the satellite diversity gain. These findings provide first analytical steps toward designing and user connection strategies to NTNs in urban environments.
\end{abstract}


\IEEEpeerreviewmaketitle
\section{Introduction}

To ensure seamless 6G coverage over terrestrial and non-terrestrial regions, NTNs such as low Earth orbit (LEO) satellites, unmanned aerial vehicles (UAVs), and high-altitude platform stations (HAPSs) are considered as essential network elements of future wireless systems. These non-terrestrial platforms can provide a complementary service to ground-based stations, especially in areas which are remote, sparsely populated, or affected by natural disasters. More generally, NTN elements are essential to ensure the vertical integration of aerial and terrestrial communication coverage \cite{geraci2022integrating}.

In dense urban environments, NTNs also have significant potential to meet the data demands of ground users and maximize communication availability. However, in such environments, line-of-sight (LoS) links between NTN elements and ground users are limited by numerous obstacles, such as high-rise buildings. The urban skyline is a key factor in determining the visibility of satellites or UAVs from the user's point of view. Therefore, it is essential to model the blocking phenomena experienced by users in a three-dimensional (3D) space and to analyze the statistical characteristics of sky visibility in this context. The aim is to then evaluate the performance and availability of NTNs in urban environments and derive methods for optimal system design \cite{wang2025modeling}.

In this paper, we propose and analyze the blocking effect by modeling the layout of buildings in urban environments using three-dimensional stochastic geometry \cite{baccelli2009stochastic}. This work explicitly extends the two-dimensional sky visibility analysis framework recently proposed \cite{lee2024much}. This previous work focused on improving the blocking angle distribution along a user's view direction using reconfigurable intelligent surfaces (RIS) under a simplified two-dimensional model. This paper generalizes this analysis of the blockage angle distribution by defining the skyline process, which continuously varies with the azimuth angle in a 3D network environment. This in particular allows one to analyze the joint blocking properties and spatial correlation across the entire angular domain.

\subsection{Related Works}
\subsubsection{Network Analysis using Stochastic Geometry}
Stochastic geometry has become a key tool in network performance analysis, providing analytical utility in dealing with the random location distributions of wireless network nodes \cite{baccelli2009stochastic}. For example, stochastic geometry has been used to characterize network performance of ad-hoc networks \cite{baccelli2006aloha,baccelli2009stochasticopp,lee2016spectral}, cellular networks \cite{andrews2011tractable,dhillon2012modeling,di2013average}, and V2X networks \cite{tong2016stochastic,yi2019modeling} in terrestrial environments. Recently, in the field of NTNs, this mathematical framework has been used to analyze network-level indicators such as visibility, connectivity, coverage, and handover rate \cite{chetlur2017downlink,banagar2020performance,okati2022coverage,al2022next}. However, these studies have mainly focused on the randomness of NTN network elements and tend to simplify the blockages caused by terrestrial environments.

\subsubsection{Urban Blockage Modeling}

Modeling blockages in terrestrial urban environments has mainly relied on empirical models, such as ITU-R models \cite{ITUR_P1411}, where the LoS probability decreases with distance between transceivers. Some studies analyze this blockage using random shape theory \cite{kendall1989survey}. By incorporating this theory, independent blockage models have been proposed to analyze the blockage effect \cite{bai2014analysis,blaszczyszyn2015wireless,ilow1998analytic}. By using a Poisson line process, correlated shadowing fields of urban networks and in-building networks have been analyzed to combine general stochastic geometry models with correlated shadowing, where links at nearby locations are blocked by the same obstacles \cite{baccelli2015correlated,zhang2015indoor,lee20163}. However, this approach is insufficient to describe the blocking distribution based on elevation that determines the availability of NTN communication links. The author in \cite{al2020probability,al2024line} analyzed LoS probability of links in an outdoor environment, and a multiplicative cascade blockage model has been incorporated to analyze the blockage effect in urban environments \cite{baccelli2022user, liu2022coverage, liu2023macro, liu20263d}.

\subsubsection{Sky Visibility and Proposed Extension}
In a recent study closely related to this paper, we proposed a framework to quantify the sky visibility experienced by users in two-dimensional urban environments using an one-dimensional marked point process \cite{lee2024much}. This study derived the blockage angle distribution in a specific direction by modeling the height information of the building as marks, demonstrating how RISs can improve visibility. However, \cite{lee2024much} has limitations in that it is mainly focused on the analysis of blocking probabilities and first-order statistics in a single direction. Especially, in real scenarios where satellites move quickly or require dual connectivity, it is essential to know the spatial correlation of blocking angles that continuously change across all azimuths, beyond the simple unidirectional analysis of \cite{lee2024much}.

Thus, this paper explicitly extends the previous model to 3D by defining the skyline process which is based on a two-dimensional marked point process. We model the blocking elevation as a stochastic process through a function $\omega(\psi)$ of azimuth $\psi$ and generalize the existing unidirectional analysis by identifying the joint distribution between any two azimuths.

\subsection{Contributions}
The main technical contributions of this paper are the following:
\begin{itemize}
\item \textbf{Mathematical Definition of Skyline Process:} By leveraging the 3D stochastic geometry framework, we propose the user-centric skyline process to analyze the ground user's sky visibility. This framework allows one to characterize the maximum blockage elevation angle as a function of the azimuth, describing the 3D sky visibility from the user's perspective.

\item \textbf{Closed-Form Distributions:} We derive closed-form expressions for the Cumulative Distribution Function (CDF) and Probability Density Function (PDF) of the maximum blockage elevation angle in any arbitrary direction $\psi$.

\item \textbf{Spatial Correlation and Joint Statistics:} We derive the joint CDF of the blockage angles for two different directions, $\psi_1$ and $\psi_2$. Furthermore, we analyze the second-order statistics of this stochastic process, including the Auto-Correlation Function (ACF) and Power Spectral Density (PSD). These metrics quantify the spatial correlation and the smoothness of the blockage process, which are essential for analyzing diversity gains and handover strategies.

\item \textbf{Application to NTN System Design:} Based on the analytical results, we analyze the performance of LEO satellite networks in urban environments. We propose an optimization framework for the minimum elevation mask angle that determines how many satellites are visible in an average sense to achieve a target outage probability, and analyze the diversity gain by leveraging the spatial correlation of the blockage effect.
\end{itemize}

\section{System Model and Geometric Formulation}\label{sec:sec2}

\subsection{Modeling of 3D Urban Environment}\label{sec:sec2-1}

\begin{figure}[t]
     \centering
     \begin{subfigure}[b]{0.42\textwidth}
         \centering
         \includegraphics[scale=0.30]{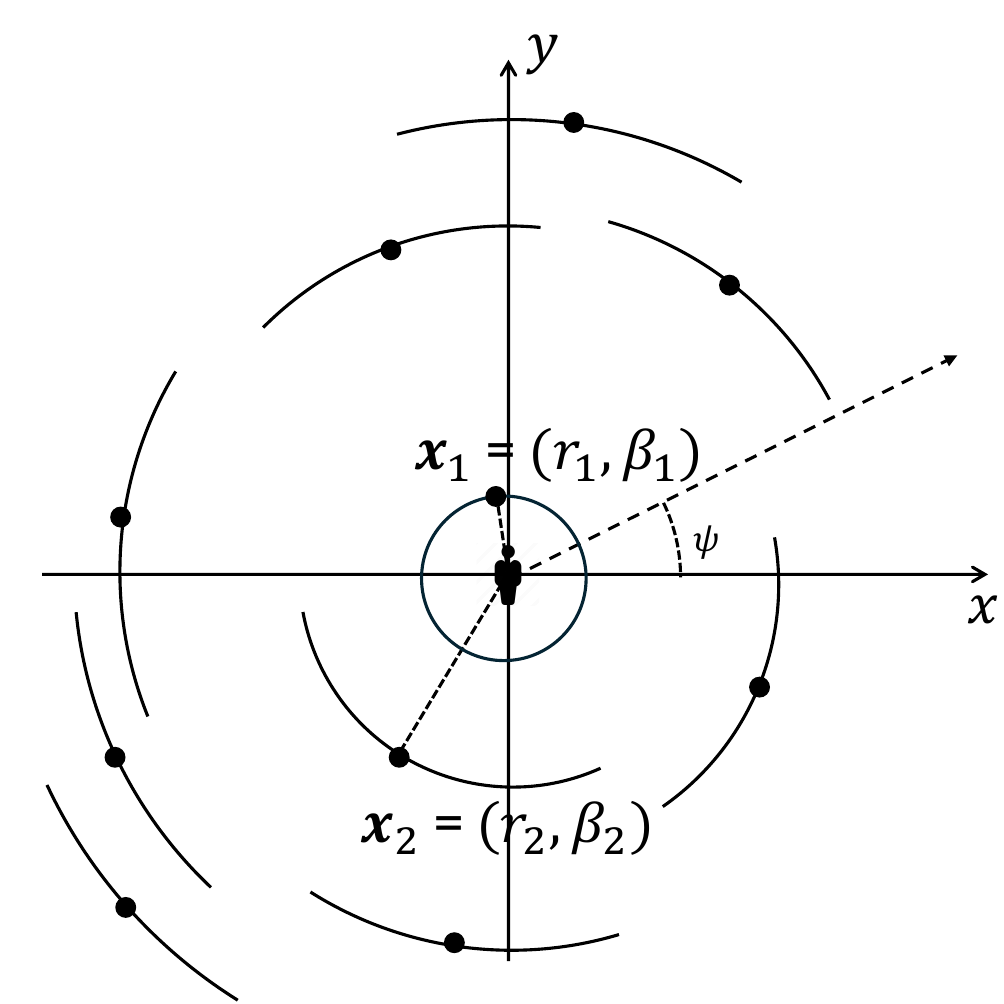}
         \caption{A top-view illustration of the 3D user-centric network model.}
         \label{fig:sys_model_3D} 
     \end{subfigure}
     \hfill
     \vspace{-1em}
     \begin{subfigure}[b]{0.42\textwidth}
         \centering
         \includegraphics[scale=0.30]{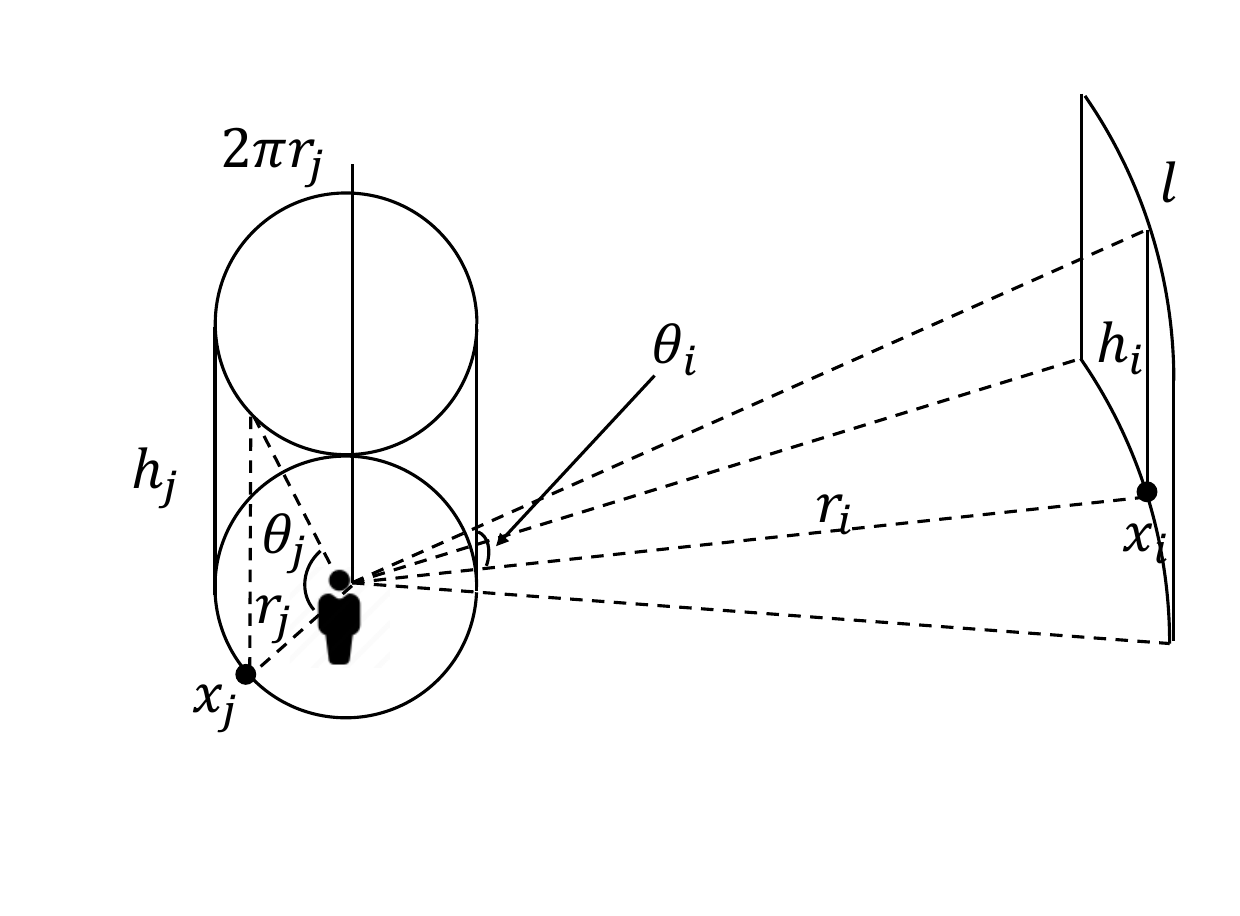}
         \caption{A parameter illustration of individual buildings.}
         \label{fig:obs_angle} 
     \end{subfigure}
     \caption{Illustration of the urban network model.} 
     \label{fig:sys_model}
\vspace{-1em}\end{figure}

In this section, we introduce the 3D urban environment model representing the spatial distribution of buildings. This model captures the LoS connectivity between a tagged user and aerial network nodes, affected by surrounding buildings. We consider a scenario where the user is located on the ground, and aim to develop a novel framework for analyzing blockage elevation angles, sky visibility and connectivity between the user and NTN nodes. This model consists of two dimensions representing the urban area (modeled as the Euclidean plane) and a third dimension representing the building heights (modeled as elements of the positive half-line). The 3D user-centric urban network model introduced here is inspired by the model defined in \cite{baccelli2022user} and is selected for mathematical convenience.\footnote{This user-centric PPP-based model introduced in \cite{baccelli2022user} is essential. The tractability comes from the fact that the upper boundary of a building is "seen" (more precisely approximated) as a constant.}

The user is located at the origin, namely $(0,0,0)$. We define the set of buildings as $B =\{b_i\}$, and each $b_i$ is modeled as an arc segment of a cylindrical surface with the center of the cylinder's base located at the origin. We further assume that the buildings $\{b_i\}$ are numbered in increasing order of their distance from the user, so that a lower index corresponds to a closer building. We denote the locations of the centers of the arc segments on the building's base as $\Phi=\{\mathbf{x}_i\}$, which are assumed to be a realization of a two-dimensional homogeneous marked Poisson point process (PPP) with intensity $\lambda$, with associated marks $H=\{h_i\}$ representing the heights. We denote by $\mu$ the inverse of the average height. The marks of $\Phi$ are assumed to be independent and identically distributed (i.i.d.) and independent of $\Phi$, with CDF $F_H(\cdot)$ and corresponding probability density function (PDF) $f_H(\cdot) = \frac{\mathrm{d}F_H(\cdot)}{\mathrm{d}h}$. We further define the parameter $\rho = \frac{\lambda}{\mu}  = \lambda \mathbb{E}[h_i]$. Since $\rho$ is the product of the building density and the average height, a higher $\rho$ represents a denser urban environment. Figure \ref{fig:sys_model_3D} illustrates a view from above of a realization of our 3D user-centric network. 

We denote by $r_i$ and $\beta_i$ the distance from the origin to $\mathbf{x}_i$ and the orientation angle relative to the $x$-axis, respectively, where $r_i\in [0,\infty)$ and $\beta_i \in (-\pi, \pi] $. In other words, $\mathbf{x}_i = (r_i \cos \beta_i, r_i \sin \beta_i) \in \mathbb{R}^2$. Further, we define $l$ as the arc length of each building.
If $r_i > \frac{l}{2\pi}$, we define the arc segment of building, $b_i$, with length $l$, centered
at the user. If $r_i \le \frac{l}{2\pi}$, the arc segment becomes the entire circle centered at the user with radius $r_i$, since its circumference
$2\pi r_i$ does not exceed $l$, as illustrated in Figure \ref{fig:obs_angle}. The $i$-th building, $b_i$, is 
\begin{align}
\begin{split}
    b_i &= \Bigg\{ (x,y,z) \bigg| x^2+y^2 = r_i^2, \mbox{Arg}(x,y)-\beta_i\in\nonumber \\
    &\qquad \left[-\min\left(\frac{l}{2r_i},\pi\right),\min\left(\frac{l}{2r_i},\pi\right)\right] ,0\leq z \leq h_i \Bigg\},
    \end{split}
\end{align}
where $\mbox{Arg}(x,y)\in(-\pi,\pi]$ is the argument of $x+jy$.

\subsection{Skyline Process: LoS Sky Visibility from the User}

In this subsection, we introduce the skyline process, which characterizes the ground user’s LoS sky visibility. Based on the distribution of buildings presented in Section~\ref{sec:sec2-1}, we first define the skyline process, which represents the upper envelope of the maximum blocking elevation angles observed by the ground user toward the tops of surrounding buildings, and then analyze its statistical properties.

\begin{figure}[!t] 
\centering
\begin{subfigure}{0.49\columnwidth}
    \centering
    \includegraphics[width=\linewidth]{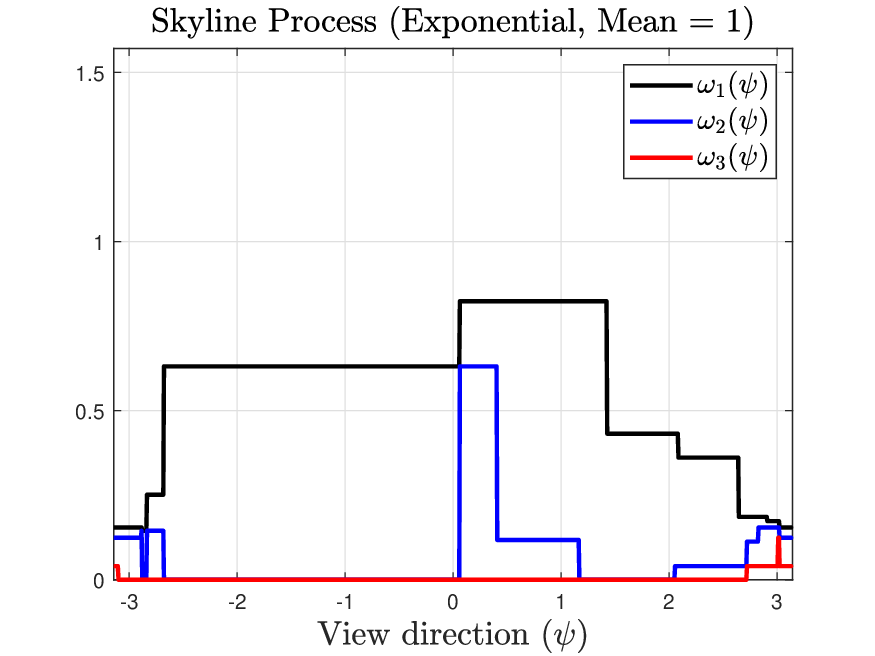}
    \caption{$\mu = 1$.}
    \label{fig:subfig1}
\end{subfigure}
\hfill
\begin{subfigure}{0.49\columnwidth}
    \centering
    \includegraphics[width=\linewidth]{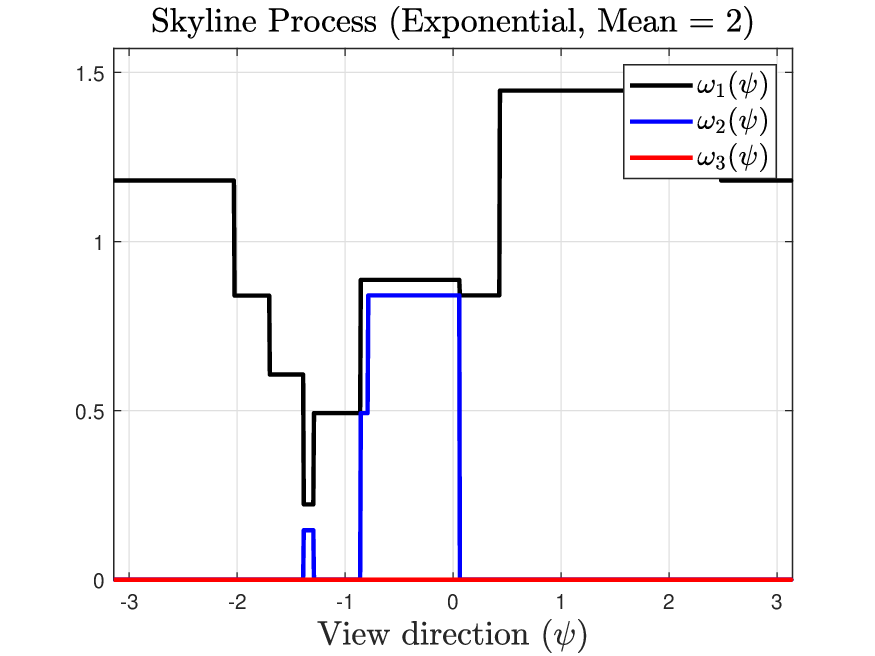}
    \caption{$\mu = 0.5$.}
    \label{fig:subfig2}
\end{subfigure}

\begin{subfigure}{0.49\columnwidth}
    \centering
    \includegraphics[width=\linewidth]{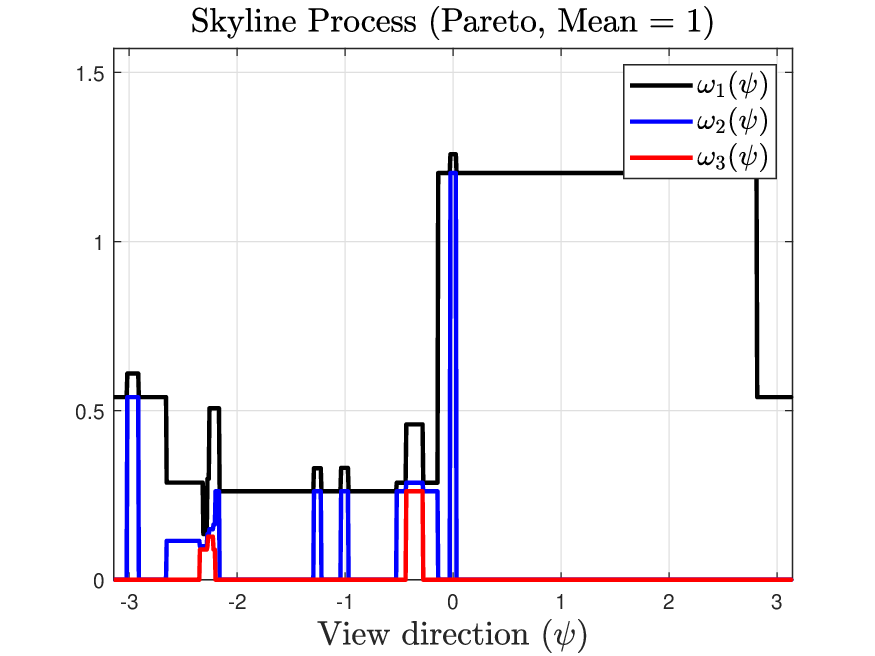}
    \caption{$\kappa = 1.5, s = \frac{1}{3}$.}
    \label{fig:subfig3}
\end{subfigure}
\hfill
\begin{subfigure}{0.49\columnwidth}
    \centering
    \includegraphics[width=\linewidth]{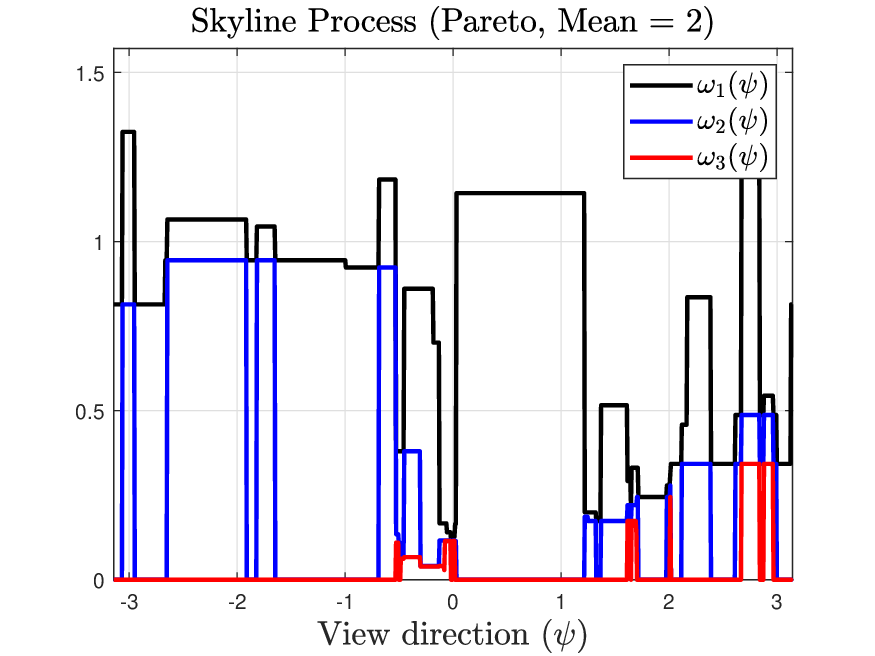}
    \caption{$\kappa = 1.5, s = \frac{2}{3}$.}
    \label{fig:subfig4}
\end{subfigure}
\vspace{-0.5em}
\caption{Illustrations of Skyline processes.}
\label{fig:slp_2x2}
\vspace{-1em}\end{figure}

Let us denote by $\theta_i$ the blockage elevation angle created by the $i$-th building, $b_i$, i.e., $\theta_i = \arctan\frac{h_i}{r_i}$. As shown in Figure \ref{fig:sys_model_3D}, $\psi\in(-\pi,\pi]$ is defined as the viewing direction of the user, which is the angle with respect to the $x$-axis. In addition, we denote by $N(\psi)$ the set of buildings that intersect the viewing direction $\psi$, and by $V(\psi)$ the subset of $N(\psi)$ consisting of buildings whose tops are visible to the user along that direction in a LoS manner. More precisely, 
\begin{align}
    N(\psi)& =  \left\{b_i\in B \bigg| |\psi-\beta_i|\le \min\left(\frac{l}{2r_i},\pi\right) \right\},\nonumber\\
    V(\psi)& = \bigg\{b_i\in N(\psi) \bigg| \theta_j < \theta_i \nonumber \\
    &\qquad \mbox{for all } j<i \mbox{ such that } b_j,b_i\in N(\psi) \bigg\}.\nonumber
\end{align}

For a given direction $\psi$, let $\omega_k(\psi)$ be the $k$-th largest blockage elevation angle among all $\{\theta_i\}$ created by buildings in $V(\psi)$. We set $\omega_k(\psi)=0$ for $k>|V(\psi)|$. Consequently, $\omega_1(\psi) = \sup_{b_i\in V(\psi)}\theta_i$, is the largest blockage elevation angle in direction $\psi$. Note that the user's LoS visibility into the sky is determined by $\omega_1(\psi)$. We define the mapping $\psi\in(-\pi,\pi] \rightarrow \omega_k(\psi)\in\left[0,\frac{\pi}{2}\right]$, for $k>0$, as the $k$-th skyline process. Figure \ref{fig:slp_2x2} gives realizations of $\omega_1(\psi)$, $\omega_2(\psi)$ and $\omega_3(\psi)$ when $\lambda = 0.1$, and $l=5$, where the building heights $h_i$ follow exponential and Pareto distributions, respectively.

\section{Statistical Properties of the Skyline Process}\label{sec:sec3}

In this section, we investigate the statistical properties of the skyline processes for the directional and omnidirectional cases. We summarize the key statistical properties derived in this paper, which mainly characterize the blockage elevation angle based visibility and blockage patterns in the 3D stochastic geometry framework.

\begin{itemize}

\item \textbf{Geometry of the Blockage Region, $\boldsymbol{\mathcal{A}(\psi)}$}, where building centroids lie to obstruct the view in direction $\psi$.

\item \textbf{CDF of $\boldsymbol{\omega_1(\psi)}$}, the maximum elevation angle in a given direction $\psi$, indicating the likelihood of LoS blockage.
(see Theorem~\ref{theo:CDF_omega1_psi}).

\item \textbf{Distribution of $\boldsymbol{r_n(\psi)}$}, the horizontal distance to the $n$-th nearest building intersecting direction $\psi$
(see Lemma~\ref{lem:n-thneighbor}).

\item \textbf{Joint Distribution of $\boldsymbol{(R(\psi), H(\psi))}$}, the distance and height of the building that creates the maximum elevation angle in direction $\psi$
(see Theorem~\ref{theo:joint_PDF_psi}).

\item \textbf{Probability that $\boldsymbol{\omega_2(\psi) = 0}$} or equivalently probability that the visibility of all buildings except the nearest building is blocked along the direction $\psi$
(see Theorem~\ref{theo:omega2}).

\item \textbf{CDF and Expectation of $\boldsymbol{\sup_{\psi} \omega_1(\psi)}$}, the global supremum of all elevation angles over all directions. The expectation is also derived in closed form for exponential heights
(see Theorem~\ref{theo:max_omega_1_psi_CDF}).

\item \textbf{Joint Distribution of $\boldsymbol{(R^*, H^*)}$}, the distance and height of the building that causes the global supremum of $\omega_1(\psi)$, across all directions
(see Theorem~\ref{theo:joint_density}).

\item \textbf{Joint Spatial Statistics for Angular Separation}, i.e., joint CDF of $\omega_1(\psi_1)$ and $\omega_1(\psi_2)$, quantifying in particular the spatial correlation of the blockage process across different viewing directions (see Theorem~\ref{theo:joint_CDF_omega_1_psi_1_psi_2}).

\item \textbf{Probability of $\boldsymbol{\omega_1(\psi_1) = \omega_1(\psi_2)}$}, i.e., likelihood that a single building is responsible for the blockage in both directions $\psi_1$ and $\psi_2$ (see Theorem~\ref{theo:prob_omega_1_psi_1_psi_2_equal}).

\item \textbf{Second-Order Statistics (ACF)}, i.e., Auto-Correlation Function, $R_{\omega_1}(\delta)$, to measure the angular dependence and smoothness of the skyline process.

\item \textbf{Spectral Analysis (PSD)}, i.e., Power Spectral Density, $S_{\omega_1}(f)$, to characterize the fluctuation rate and roughness of the blockage angle in the spatial frequency domain.
\end{itemize}%
\vspace{-1.0em}

\subsection{First-Order Statistics: Blockage Probability} 
\subsubsection{Skyline Process in a Given Direction}
First, we investigate the distribution of the buildings that intersect the user’s view direction $\psi \in (-\pi, \pi]$, relative to the $x$-axis to analyze the statistical properties of the skyline process.  

\begin{figure}[t]
\begin{center}
\includegraphics[width=0.78\columnwidth]{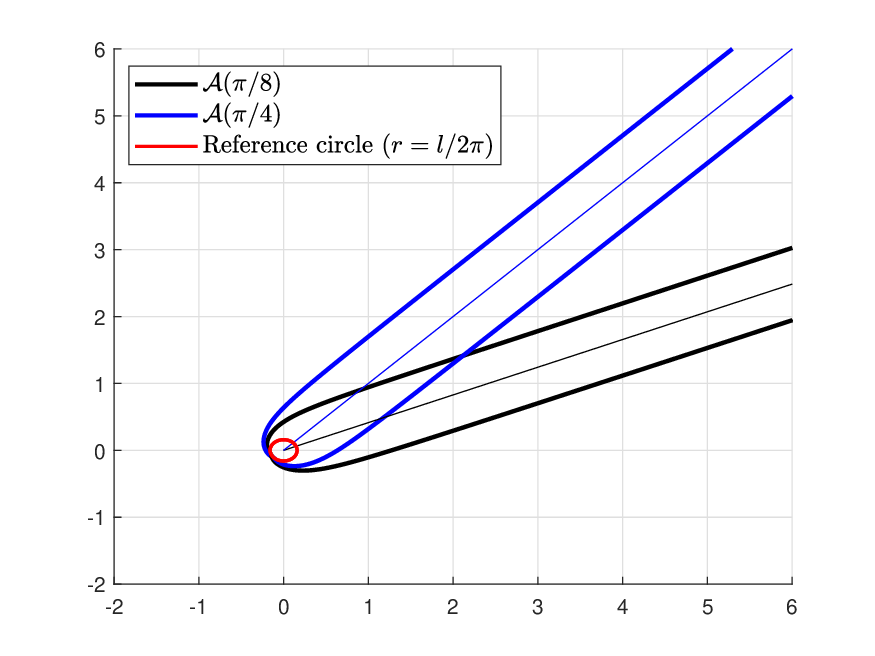}
\end{center}
\vspace{-0.5em} 
\caption{Illustration of $\mathcal{A}(\psi)$ where $\psi = \frac{\pi}{8}, \frac{\pi}{4}$ and $l=1$.} 
\label{fig:region_crossing}
\vspace{-1.2em}
\end{figure}

\noindent \textbf{Geometry of the Blockage Region:}
A building whose center has polar coordinates $(r, \alpha)$, i.e., $\mathbf{x}_i=(r \cos \alpha, r\sin\alpha)$, intersects the user’s view direction $\psi$ if and only if $\mathbf{x}_i$ is in 
\begin{align}
    \mathcal{A}(\psi) &= \left\{ (r, \alpha)  \bigg|  r |\alpha - \psi| \leq \frac{l}{2},   r > 0 \right\}.\nonumber
\end{align}
This is equivalent to
\begin{small}
\begin{align}
    \mathcal{A}(\psi) = 
    \begin{cases}
        \left\{ (r, \alpha) \big| \alpha \in [-\pi, \pi] \right\} & \text{for } r < \frac{l}{2\pi}, \\
        \bigg\{ (r, \alpha) \bigg| \alpha \in \left[ \psi - \frac{l}{2r}, \psi + \frac{l}{2r} \right] \\
        \qquad \mod 2\pi \bigg\} & \text{for } r \geq \frac{l}{2\pi}.
    \end{cases}\nonumber
\end{align}
\end{small}
Figure \ref{fig:region_crossing} depicts $\mathcal{A}(\psi)$ when $\psi = \frac{\pi}{8}, \frac{\pi}{4}$ and $l=1$. The following theorem gives the distribution of $\mathbf{x}_i$ within $\mathcal{A}(\psi)$.

\theorem The number of buildings intersecting the user's viewing direction $\psi$ and within a distance $R$ from the user follows a Poisson distribution with parameter $\lambda l\left(R - \frac{l}{4\pi}\right)$ for $R > \frac{l}{2\pi}$ and $\lambda \pi R^2$ otherwise.
\begin{IEEEproof}
    Let $B(0,R)$ denote the ball of radius $R$ centered at the user. A building intersects the user's view direction $\psi$ if its centroid lies in the region $\mathcal{A}(\psi)\cap B(0,R)$. The area of $\mathcal{A}(\psi)\cap B(0,R)$ is 
\begin{align} 
    &|\mathcal{A}(\psi)\cap B(0,R)| \nonumber \\
    &= \begin{cases}
        \begin{aligned}
            &\int_{0}^{\frac{l}{2\pi}}\int_{-\pi}^{\pi} r\,\mathrm{d}\theta\,\mathrm{d}r + \int_{\frac{l}{2\pi}}^{R}\int_{\psi-\frac{l}{2r}}^{\psi+\frac{l}{2r}} r\,\mathrm{d}\theta\,\mathrm{d}r \\
            &\quad = l\left(R-\frac{l}{4\pi}\right)
        \end{aligned} & \text{for } R > \frac{l}{2\pi}, \\
        \int_{0}^{R}\int_{-\pi}^{\pi} r\,\mathrm{d}\theta\,\mathrm{d}r = \pi R^2 & \text{for } R \leq \frac{l}{2\pi}.
    \end{cases}\nonumber
\end{align}%
So, the number of buildings intersecting the user's view direction $\psi$ within a distance $R$, follows a Poisson distribution with parameter $\lambda l\left(R - \frac{l}{4\pi}\right)$ for $R > \frac{l}{2\pi}$ and $\lambda \pi R^2$ for $R \leq \frac{l}{2\pi}$.
\end{IEEEproof}

\noindent\textbf{CDF of $\boldsymbol{\omega_1(\psi)}$:} Here, $\omega_1(\psi)$ denotes the maximum elevation blockage angle seen by the user in the direction $\psi$, and is defined as the arctangent of the ratio of height and distance of the most building blocking the LoS visibility toward non-terrestrial nodes along the viewing direction.

\theorem\label{theo:CDF_omega1_psi}
For a given direction $\psi \in (-\pi,\pi]$, the CDF of $\omega_1(\psi)$ is
\begin{align} \label{eq:CDF_omega_1}
    &\mathbb{P}\bigl[\omega_1(\psi) \leq \phi\bigr]= \exp\Biggl(
        - \lambda \biggl(
        2\pi \int_{0}^{\frac{l}{2\pi}} \bigl(1 - F_H(r \tan \phi)\bigr) \, r \, \mathrm{d}r \nonumber \\
        &\qquad + \int_{\frac{l}{2\pi}}^{\infty}
          \int_{-\frac{l}{2r} + \psi}^{\frac{l}{2r} + \psi}
          \bigl(1 - F_H(r \tan \phi)\bigr) \, r \, \mathrm{d}\theta \, \mathrm{d}r
        \biggr)
    \Biggr).
\end{align}

\begin{IEEEproof}
    For a given $\psi$, 
    \begin{align}
    &\mathbb{P}[\tan\omega_1(\psi)\leq\tau]=\mathbb{P}\left[\sup_{\mathbf{x}_i\in\Phi\cap\mathcal{A}(\psi)}\frac{h_i}{r_i}\leq \tau\right]\nonumber\\
    &=\mathbb{E}\left[\mathbbm{1}\left(\sup_{\mathbf{x}_i\in\Phi\cap\mathcal{A}(\psi)}\frac{h_i}{r_i}\leq \tau\right)\right]=\mathbb{E}\left[\prod_{\mathbf{x}_i\in\Phi\cap\mathcal{A}(\psi)}\mathbbm{1}\left(\frac{h_i}{r_i}\leq \tau\right)\right]\nonumber\\        
    &=\mathbb{E}\left[\exp\left(\sum_{\mathbf{x}_i\in\Phi\cap\mathcal{A}(\psi)}\log\left(\mathbbm{1}\left(\frac{h_i}{r_i}\leq \tau\right)\right)\right)\right]\nonumber\\
    &\stackrel{(a)}{=}\exp\Biggl(-\lambda\biggl(\int_{0}^{\frac{l}{2\pi}}\int_{-\pi}^{\pi}\left(1-\mathbb{E}_h\left[\mathbbm{1}\left(\frac{h}{r}<\tau\right)\right]\right) r\,\mathrm{d}\theta\,\mathrm{d}r\nonumber\\
    &\qquad+\int_{\frac{l}{2\pi}}^{\infty}\int_{-\frac{l}{2r}+\psi}^{\frac{l}{2r}+\psi}\left(1-\mathbb{E}_h\left[\mathbbm{1}\left(\frac{h}{r}<\tau\right)\right]\right) r\,\mathrm{d}\theta\,\mathrm{d}r\biggr)\Biggr)\nonumber\\
    &=\exp\Biggl(-\lambda\biggl(2\pi\int_{0}^{\frac{l}{2\pi}}\bigl(1-F_H(r\tau)\bigr) r\,\mathrm{d}r\nonumber\\
    &\qquad +\int_{\frac{l}{2\pi}}^{\infty}\int_{-\frac{l}{2r}+\psi}^{\frac{l}{2r}+\psi}\bigl(1-F_H(r\tau)\bigr) r\,\mathrm{d}\theta\,\mathrm{d}r\biggr)\Biggr),\nonumber
\end{align}
    where (a) comes from the probability generating functional (PGFL) of the PPP \cite{baccelli2009stochastic}. Since $\mathbb{P}[\tan\omega_1(\psi)\leq \tau] = \mathbb{P}[\omega_1(\psi)\leq \arctan \tau]$, we obtain \eqref{eq:CDF_omega_1} by substituting $\tan \phi$ for $\tau$. 
\end{IEEEproof}

\example Let $F^{\mbox{exp}}_{\omega_1}(\phi)$ denote the CDF $\omega_1(\psi)$ under the exponential distribution with parameter $\mu$. Then, \eqref{eq:CDF_omega_1} reduces to 
\begin{align}\label{eq:ex_CDF_omega1_exp}
    F^{\mbox{exp}}_{\omega_1}(\phi)=e^{-\frac{2\pi\lambda\left(1-\exp\left(-\frac{l\mu\tan\phi}{2\pi}\right)\right)}{\mu^2(\tan\phi)^2}}.
\end{align}

\example Let $F^{\mbox{par}}_{\omega_1}(\phi)$ denote the CDF $\omega_1(\psi)$ under the Pareto distribution with scale parameter $s >0$ and shape parameter $\kappa$ where $1<\kappa<2$. In this case, $F_H(h) = 1-\left(\frac{s}{h}\right)^{\kappa}$ for $h \ge s$. Then, \eqref{eq:CDF_omega_1} reduces to 
\begin{align}\label{eq:ex_CDF_omega1_par}
    F^{\mbox{par}}_{\omega_1}(\phi)=e^{\frac{\lambda l^2 (2\pi)^{-1+\kappa} \left(\frac{s}{ l\tan{\phi}}\right)^\kappa}{\kappa^2 - 3\kappa + 2}
}.
\end{align}

\begin{figure}[t]
\begin{center}
\includegraphics[width=0.9\columnwidth]{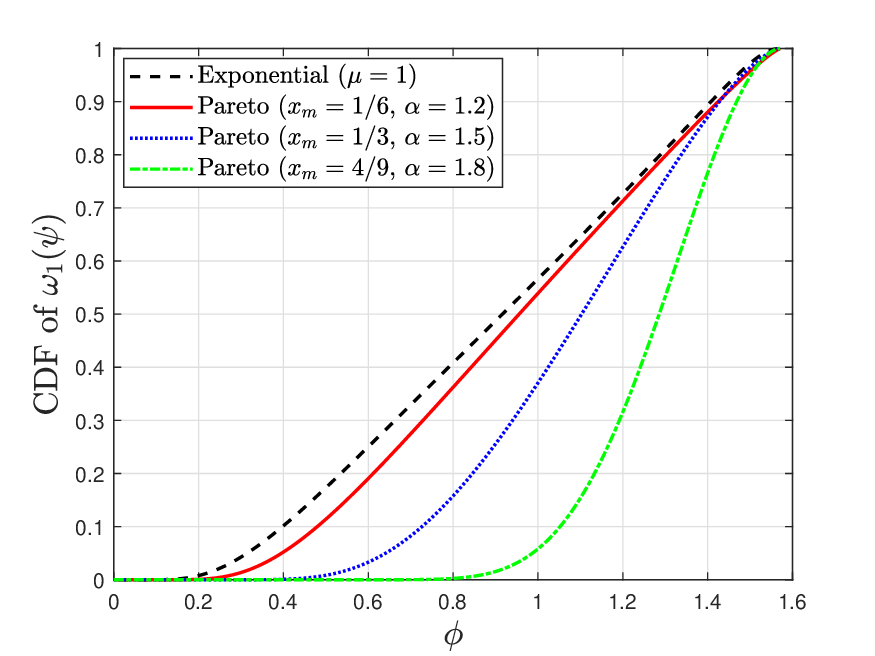}
\end{center}
\vspace{-0.5em} 
\caption{CDF of $\omega_1(\psi)$ with $\lambda = 1$ and $l = 1$.} 
\label{fig:CDF_omega_1}
\vspace{-1em}\end{figure}

Figure \ref{fig:CDF_omega_1} illustrates the CDF of $\omega_1(\psi)$ under exponential and Pareto height distributions with $\lambda = 1$ and $l = 1$. In this figure, we set the same average building height to $1$ with $\mu=1$ and $s=\frac{\kappa-1}{\kappa}$. In this figure, we can observe the relation $F^{\mbox{exp}}_{\omega_1}(\phi)>F^{\mbox{par}}_{\omega_1}(\phi)$, when the mean building heights are the identical. This indicates a stochastic domination of the blockage angles, reflecting the impact of the heavy-tailed building height distribution. The proof of $F^{\mbox{exp}}_{\omega_1}(\phi)>F^{\mbox{par}}_{\omega_1}(\phi)$ for all $\phi$, is provided in Appendix~\ref{appen:eq1213}. 

\noindent\textbf{Distribution of $\boldsymbol{r_n(\psi)}$:}
\lemma\label{lem:n-thneighbor} Let $r_n(\psi)$ denote the distance from the user to the $n$-th nearest building (visible or not) that intersects the user's viewing direction $\psi$. Then, the PDF of $r_n(\psi)$ is
\begin{align}\label{eq:ndist_Apsi}
    &f_{r_n(\psi)}(r) \nonumber \\
    &= \begin{cases}
        2\lambda\pi r\frac{(\lambda\pi r^2)^{n-1}}{(n-1)!}e^{-\lambda\pi r^2}& \text{for } r \leq \frac{l}{2\pi}, \\[2.5ex]
        \frac{\lambda l \left(\lambda l \left(r-\frac{l}{4\pi}\right)\right)^{n-1}}{(n-1)!}e^{-\lambda l \left(r-\frac{l}{4\pi}\right)} & \text{for } r > \frac{l}{2\pi}.
    \end{cases}
\end{align}
\begin{IEEEproof}
    Since $\{\mathbf{x}_i\}$ is a realization of a homogeneous PPP with intensity $\lambda$, the number of $\mathbf{x}_i$ in $\mathcal{A}(\psi) \cap B(0, R)$ follows a Poisson distribution with mean $\lambda |\mathcal{A}(\psi) \cap B(0, R)|$.

The CDF of $r_n(\psi)$ is given by the probability that there are more than $n-1$ points in the region:
\begin{align}
&\mathbb{P}[r_n(\psi) \le R]  \nonumber\\
&= 1 -\sum_{k=0}^{n-1} \frac{\left( \lambda |\mathcal{A}(\psi) \cap B(0, R)| \right)^k}{k!} e^{-\lambda |\mathcal{A}(\psi) \cap B(0, R)|}.\nonumber
\end{align}
    We obtain \eqref{eq:ndist_Apsi} by taking the derivative of this CDF with respect to $R$. 
\end{IEEEproof}

Due to the rotational invariance of the building distribution with respect to the user's view direction, the distribution of $r_n(\psi)$ is independent of $\psi$.

\example By setting $n=1$, the PDF of $r_1(\psi)$ becomes
\begin{align}
    f_{r_1(\psi)}(r) = 
    \begin{cases}
        2\lambda\pi r \exp(-\lambda\pi r^2)& \mbox{for } r\leq \frac{l}{2\pi},\\
        \lambda l \exp\left(-\lambda l \left(r-\frac{l}{4\pi}\right)\right)& \mbox{for } r >\frac{l}{2\pi}.
    \end{cases}\nonumber
\end{align}

\noindent\textbf{Joint Distribution of $\boldsymbol{(R(\psi), H(\psi))}$:} Let us consider the building that creates the maximum elevation angle $\omega_1(\psi)$ in the direction $\psi$. Let $R(\psi)$ and $H(\psi)$ be the horizontal distance from the user to that building and the height of the building, respectively. The following theorem gives the joint distribution of $(R(\psi),H(\psi))$.

\theorem\label{theo:joint_PDF_psi} For any $\psi \in [-\pi, \pi)$, the joint probability density function (PDF) of the distance $R(\psi)$ and height $H(\psi)$ of the building that determines $\omega_1(\psi)$ is given by
\begin{align}
    &\mathbb{P}[R(\psi) \in [r, r + \mathrm{d}r], H(\psi) \in [h, h + \mathrm{d}h]] \nonumber \\
    &= \begin{cases}
        \begin{aligned}
            &2 \lambda \pi r \exp \biggl( -\lambda \int_{\mathcal{A}(\psi)} \Bigl[ 1 - F_H\Bigl( \frac{h \|\mathbf{x}\|}{r} \Bigr) \Bigr] \mathrm{d}\mathbf{x} \biggr) \\
            &\quad \times f_H(h) \,\mathrm{d}r \,\mathrm{d}h, 
        \end{aligned} & r \leq \frac{l}{2\pi}, \\[3ex]
        \begin{aligned}
            &\lambda l \exp \biggl( -\lambda \int_{\mathcal{A}(\psi)} \Bigl[ 1 - F_H\Bigl( \frac{h \|\mathbf{x}\|}{r} \Bigr) \Bigr] \mathrm{d}\mathbf{x} \biggr) \\
            &\quad \times f_H(h) \,\mathrm{d}r \,\mathrm{d}h, 
        \end{aligned} & r > \frac{l}{2\pi}.
    \end{cases}\nonumber
\end{align}

\begin{IEEEproof} 
{\color{black}In order to obtain the joint PDF of $R(\psi)$ and $H(\psi)$, we consider the joint Laplace transform defined as $\mathcal{L}(s,t) = \mathbb{E}[e^{-sR(\psi) - tH(\psi)}]$, which is
\begin{equation} \label{eq:jLaplace}
\begin{aligned}
    &\mathcal{L}(s,t) = \mathbb{E}\left[ \sum\nolimits_{(\mathbf{x}_i, h_i) \in \Phi} e^{-s\|\mathbf{x}_i\| - t h_i} \mathbbm{1}\left( \tfrac{h_j}{\|\mathbf{x}_j\|} \le \tfrac{h_i}{\|\mathbf{x}_i\|}, \forall j \neq i \right) \right] \\
    &\stackrel{(a)}{=} \int_{\mathcal{A}(\psi)} \int_{0}^{\infty} e^{-s\|\mathbf{x}\| - t h} \mathbb{P}^{\mathbf{x}, h} \left[ \text{No points satisfy } \tfrac{h'}{\|\mathbf{y}\|} > \tfrac{h}{\|\mathbf{x}\|} \right]\nonumber\\
    &\times f_H(h) \Lambda(\mathrm{d}\mathbf{x}) \mathrm{d}h.
\end{aligned}
\end{equation}
where (a) follows from Campbell's theorem for marked PPPs \cite{baccelli2009stochastic}, and $\mathbb{P}^{\mathbf{x}, h}$ denotes the Palm probability given a point at $(\mathbf{x}, h)$. According to Slivnyak's theorem \cite{baccelli2009stochastic}, the reduced Palm distribution is identical to the original distribution. 
So, from \eqref{eq:jLaplace}, we can identify the joint PDF $f_{R, H}(r, h)$ as the product of the spatial intensity and the void probability.

First, the spatial intensity term, obtained by integrating $\Lambda(\mathrm{d}\mathbf{x})$ over the angular component at distance $r$, is denoted as $\Lambda_{\text{loc}}(r)$, which is
\begin{align}
    \Lambda_{\text{loc}}(r) = 
    \begin{cases}
        2\pi r \lambda, & r \leq \frac{l}{2\pi}, \\
        \lambda l, & r > \frac{l}{2\pi}.
    \end{cases}\nonumber
\end{align}

Second, we calculate the probability that the region
\begin{align}
    \mathcal{S}_{r,h} = \left\{ (\mathbf{y}, h') \in \mathcal{A}(\psi) \times \mathbb{R}^+ : \frac{h'}{\|\mathbf{y}\|} > \frac{h}{r} \right\},\nonumber
\end{align}
is empty. The average number of buildings in this region is
\begin{align}
    \Lambda(\mathcal{S}_{r,h})& = \int_{\mathcal{A}(\psi)} \int_{\frac{h\|\mathbf{y}\|}{r}}^{\infty} \lambda f_H(h') \mathrm{d}h' \mathrm{d}\mathbf{y}\nonumber\\
    &= \lambda \int_{\mathcal{A}(\psi)} \left[ 1 - F_H\left( \frac{h \|\mathbf{y}\|}{r} \right) \right] \mathrm{d}\mathbf{y}.\nonumber
\end{align}
Since the number of points follows a Poisson distribution, the probability that $\mathcal{S}_{r,h}$ contains no points is given by
\begin{align}
    \mathbb{P}[\text{No points in } \mathcal{S}_{r,h}] = \exp\left( - \Lambda(\mathcal{S}_{r,h}) \right).\nonumber
\end{align}

Finally, the joint PDF is
\begin{align}
    f_{R, H}(r, h) = \Lambda_{\text{loc}}(r) f_H(h) \exp\left( - \Lambda(\mathcal{S}_{r,h}) \right),\nonumber
\end{align}
which completes the proof.}
\end{IEEEproof}

\noindent\textbf{Probability that $\boldsymbol{\omega_2(\psi) = 0}$:}
By definition, $\omega_1(\psi)$ is generated by the building which creates the maximum elevation blockage angle along the direction $\psi$, and it is always positive for any direction $\psi\in (-\pi,\pi]$. However, $\omega_2(\psi)$ can be zero when no other buildings behind the nearest building are visible. 

\theorem\label{theo:omega2} The probability that $\omega_2(\psi)=0$, for a given $\psi\in(-\pi,\pi]$, is 
\begin{align}\label{eq:omega_2_psi_zero_1}
    &\mathbb{P}[\omega_2(\psi)=0] 
    = \int_{0}^{\infty} f_H(h) \Biggl[ \int_{0}^{\frac{l}{2\pi}} 2\lambda\pi r e^{-\lambda\pi r^2}  \nonumber \\
    &\times \mathcal{V}(r, h) \,\mathrm{d}r+ \int_{\frac{l}{2\pi}}^{\infty} \lambda l e^{-\lambda l\Bigl(r-\frac{l}{4\pi}\Bigr)} \mathcal{V}(r, h) \,\mathrm{d}r \Biggr] \mathrm{d}h,
\end{align}
where $\mathcal{V}(r, h) = \exp\left(-\lambda\int_{\mathcal{A}(\psi)\setminus B(0,r)} \left[1- F_H\left(\frac{h}{r} \|\mathbf{u}\| \right)\right] d\mathbf{u} \right)$ denotes the void probability given the nearest building at distance $r$ with height $h$.

\begin{IEEEproof}
    Let $\mathbf{x}_1(\psi)$ and $h_1$ denote the {location} and the height, respectively, of the building nearest to the user in the direction $\psi$. Let $r_1 = \|\mathbf{x}_1(\psi)\|$ be the distance to this building. For given $r_1$ and $h_1$, the probability that the nearest building blocks all buildings located farther than $r_1$ is equivalent to the event that no building exists in the region $\mathcal{A}(\psi)\setminus B(0,r_1)$ with an elevation angle larger than $h_1/r_1$. This probability is given by:
\begin{align}\label{eq:tmp11_1}
    &\mathbb{P}\Biggl[\sup_{\mathbf{x}_i\in\Phi\cap(\mathcal{A}(\psi)\setminus B(0,r_1) )}\frac{h_i}{\|\mathbf{x}_i\|} < \frac{h_1}{r_1} \Biggm| r_1, h_1 \Biggr] \nonumber \\
    &= \exp\Biggl(-\lambda\int_{\mathcal{A}(\psi)\setminus B(0,r_1)} \biggl[ 1- F_H\biggl(\frac{h_1}{r_1} \|\mathbf{u}\| \biggr) \biggr] \,\mathrm{d}\mathbf{u} \Biggr).
\end{align}
    We assume that the distance to the nearest building $r_1$ follows the distribution derived from the specific point process geometry (separated into two regimes: $0 < r_1 \le l/2\pi$ and $r_1 > l/2\pi$). Finally, \eqref{eq:omega_2_psi_zero_1} is obtained by integrating \eqref{eq:tmp11_1} over the distribution of $r_1$ and $h_1$.
\end{IEEEproof}

\begin{figure}[t!]
    \centering
    \includegraphics[width=0.78\linewidth]{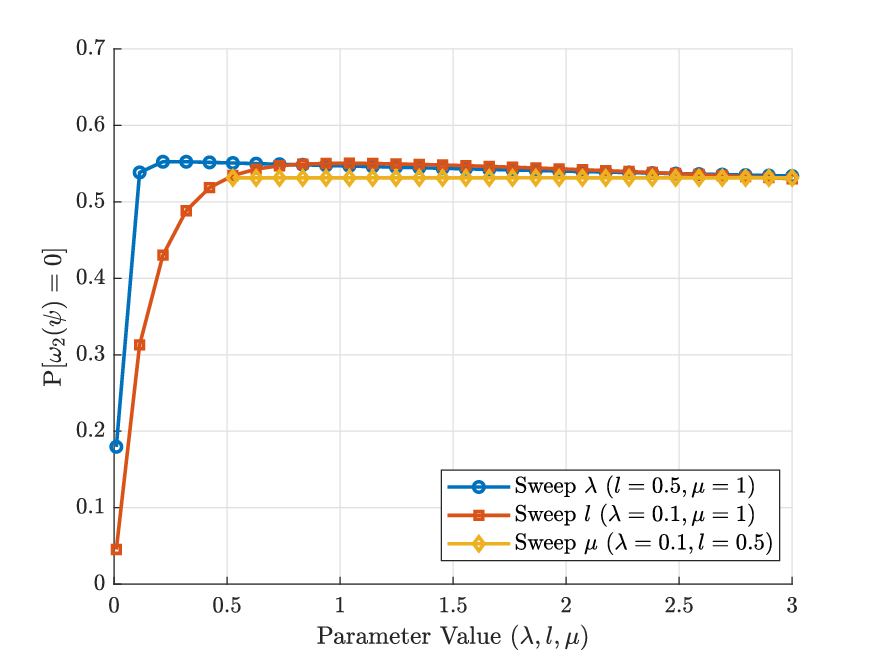}

    \caption{The probability $\mathbb{P}[\omega_2(\psi)=0]$ as a function of various urban parameters. }
    
    \label{fig:theo4}
\vspace{-1em}\end{figure}

Theorem~\ref{theo:omega2} shows that $\mathbb{P}[\omega_2(\psi) = 0]$ depends not only on the statistical properties of building heights, but also on geometric spatial parameters. Figure~\ref{fig:theo4} illustrates this probability for an exponential building height distribution with mean $1/\mu$, examining its variation with respect to the building density $\lambda$, the arc-length parameter $l$, and $\mu$.

First, when increasing the density $\lambda$, the probability increases and eventually saturates. This is because the higher the density, the greater the likelihood that the building closest to the user will be located very close, and this proximity creates a steep blockage elevation angle $\omega_1(\psi)$. 

Second, for the arc-length parameter $l$, the probability is negligible when $l$ is small but rises rapidly as $l$ increases, showing a similar saturation trend. A larger $l$ implies a wider building cross-section, which naturally casts a broader shadow, reducing the chance of observing a second independent peak.

Finally, changing the rate parameter $\mu$ has no effect on this probability. This is due to the scale-invariance property of the exponential distribution, which means that the blocking is determined by the ratio of the heights, so the absolute scale of the height of the building (all the buildings getting bigger or smaller together) cancels out.

\subsubsection{Skyline Process over All Directions}

\noindent\textbf{\\CDF and Expectation of $\boldsymbol{\sup_{\psi} \omega_1(\psi)}$:}
The following theorem provides the distribution of the maximum of $\omega_1(\psi)$. 
\theorem\label{theo:max_omega_1_psi_CDF} The CDF of the maximum of $\omega_1(\psi)$ over all azimuth angles $\psi\in(-\pi,\pi]$, is given by
\begin{align}
    &\mathbb{P}\biggl[\sup_{\psi\in(-\pi,\pi]}\omega_1(\psi)\leq \phi\biggr] =e^{-2\pi\lambda\int_{0}^{\infty} \bigl(1-F_H(r\tan\phi)\bigr) r \,\mathrm{d}r}.\nonumber
\end{align}
\begin{IEEEproof}
    The CDF of the tangent of $\sup_{\psi\in(-\pi,\pi]}\omega_1(\psi)$ is
    \begin{align}
        &\mathbb{P}\left[\tan \sup_{\psi\in(-\pi,\pi]}\omega_1(\psi)\leq\tau\right]\nonumber\\
        &=\mathbb{P}\left[\sup_{\psi\in(-\pi,\pi]}\tan \omega_1(\psi)\leq\tau\right]=\mathbb{P}\left[\sup_{\mathbf{x}_i\in\Phi}\frac{h_i}{r_i}\leq \tau\right]\nonumber\\
        &=\mathbb{E}\left[\mathbbm{1}\left(\sup_{\mathbf{x}_i\in\Phi}\frac{h_i}{r_i}\leq \tau\right)\right]=\mathbb{E}\left[\prod_{\mathbf{x}_i\in\Phi}\mathbbm{1}\left(\frac{h_i}{r_i}\leq \tau\right)\right]\nonumber\\    
        &=\exp\left(-\lambda\left(\int_{0}^{\infty}\int_{-\pi}^{\pi}\left(1-\mathbb{E}_h\left[\mathbbm{1}\left(\frac{h}{r}<\tau\right)\right]\right)            r\mbox{d}\theta \mbox{d}r\right)\right)\nonumber\\
        &=\exp\left(-2\pi\lambda\int_{0}^{\infty} (1-F_H(r\tau))r\mbox{d}r \right).\nonumber
    \end{align}
    Since $\mathbb{P}[\tan\sup_{\psi\in(-\pi,\pi]}\omega_1(\psi)\leq \tau]=$  $\mathbb{P}[\sup_{\psi\in(-\pi,\pi]}\omega_1(\psi)\leq \arctan \tau]$, the CDF of $\sup_{\psi\in(-\pi,\pi]}\omega_1(\psi)$ becomes
    \begin{align}
    e^{-2\pi\lambda\int_{0}^{\infty} (1-F_H(r\tan \phi))r\mbox{d}r},\nonumber
    \end{align}
    by substituting $\tan \phi$ for $\tau$. 
\end{IEEEproof}

\example\label{ex:ex_CDF_max_omega1_exp} Let $\{h_i\}$ follow an exponential distribution with parameter $\mu$. Then, for $\phi \in [0,\frac{\pi}{2})$, 
\begin{align}\label{eq:ex_CDF_max_omega1_exp}
    \mathbb{P}\left[\sup_{\psi\in(-\pi,\pi]}\omega_1(\psi)\leq \phi\right] =\exp\left(-\frac{2\pi\lambda}{\mu^2 (\tan\phi)^2}\right),
\end{align}
and its PDF is
\begin{align}
    f_{\sup_{\psi\in(-\pi,\pi]}\omega_1(\psi)}(\phi) =e^{-\frac{2\pi\lambda}{\mu^2 (\tan\phi)^2}}\frac{4\pi\lambda}{\mu^2\tan \phi (\sin \phi)^2}.\nonumber
\end{align}

There is an obvious stochastic dominance between $\omega_1(\psi)$ and $\sup_{\psi\in(-\pi,\pi]}\omega_1(\psi)$. Equation \eqref{eq:ex_CDF_omega1_exp} converges to \eqref{eq:ex_CDF_max_omega1_exp} as $l \to \infty$. We can observe that the distribution of the maximum blockage elevation angle over all $\psi$ is equivalent to that for a given $\psi$ in a scenario where the length of buildings goes to infinity. The discrepancy between the two equations comes from the visibility gain arising from the gaps created by the finite width.

\example In Example \ref{ex:ex_CDF_max_omega1_exp}, 
\begin{align}\label{eq:mean_max_omega}
    \mathbb{E}\left[\sup_{\psi\in(-\pi,\pi]}\omega_1(\psi)\right] = \frac{\pi}{2}-\frac{\pi}{2}e^{\frac{2\pi\lambda}{\mu^2}}\mbox{Erfc}\left(\frac{\sqrt{2\pi\lambda}}{\mu}\right),
\end{align}
where $\mbox{Erfc}(\cdot)$ is the complementary error function  
\begin{align}
    \mbox{Erfc}(z) = 1-\frac{2}{\sqrt{\pi}}\int_{0}^{z}e^{-t^2}\mathrm{d}t.\nonumber
\end{align}

\begin{figure}[t] 
    \centering
    
    \begin{subfigure}{1.0\columnwidth}
        \centering
        \includegraphics[width=0.78\linewidth]{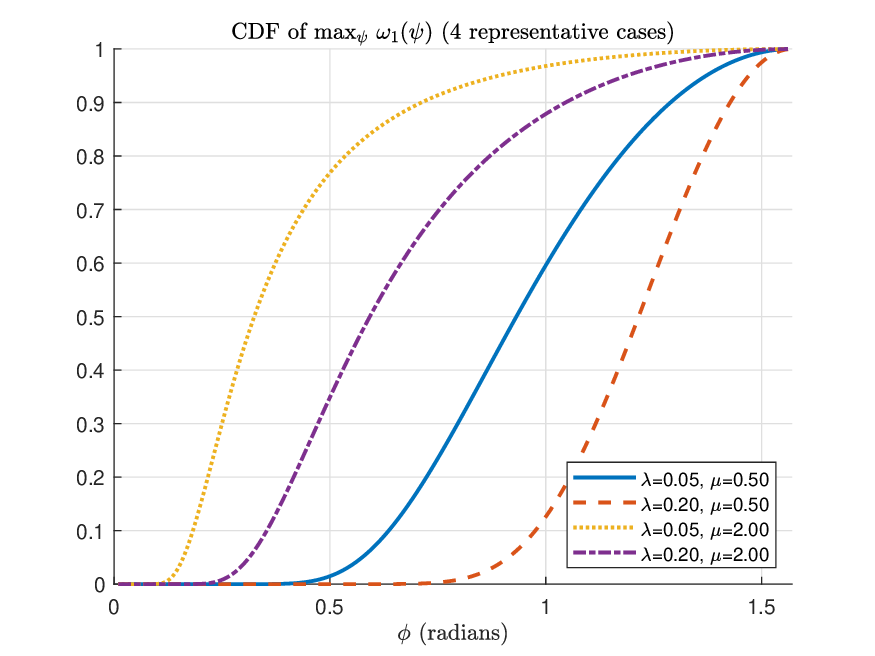}
        \caption{CDF of $\sup_{\psi}\omega_1(\psi)$}
        \label{fig:eq26}
    \end{subfigure}
\vspace{-0.5em}

    \begin{subfigure}{1.0\columnwidth}
        \centering
        \includegraphics[width=0.78\linewidth]{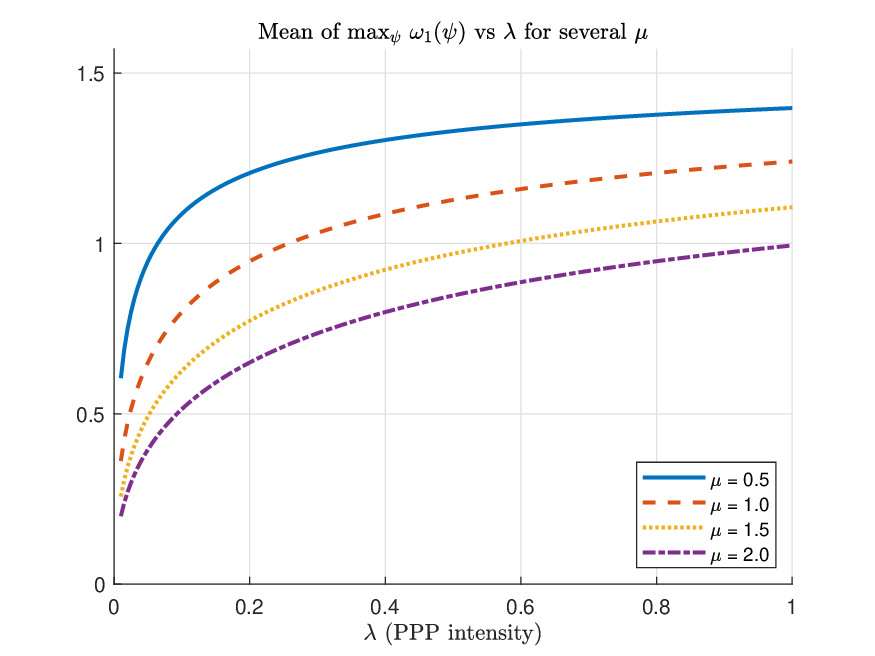}
        \caption{Mean of $\sup_{\psi}\omega_1(\psi)$}
        \label{fig:eq28}
    \end{subfigure}

    \caption{CDF and Mean of $\sup_{\psi}\omega_1(\psi)$}
    \label{fig:eq26eq28}
\vspace{-1em}\end{figure}

Figure~\ref{fig:eq26eq28} illustrates the CDF and mean of the maximum elevation angle, $\sup_{\psi}\omega_1(\psi)$. As shown in Fig.~\ref{fig:eq26}, the CDF curve shifts to the right as $\lambda$ increases or $\mu$ decreases. This implies that users in denser or taller urban environments are statistically more likely to experience steeper blockage angles. Fig.~\ref{fig:eq28} depicts the mean of $\sup_{\psi}\omega_1(\psi)$. We observe that the mean angle increases rapidly in the low-density regime but saturates as the buildings become denser. Notably, the impact of building height is more significant than that of density, which is consistent with the analytical result in \eqref{eq:mean_max_omega}. 

\subsection{Characterization of the Dominant Obstacle} \label{subsec:3-2}
\noindent\textbf{Joint Distribution of $\boldsymbol{(R^*, H^*)}$:}
We now focus on the building that creates the global maximum of $\omega_1(\psi)$ for $\psi\in(-\pi,\pi]$. 

\theorem\label{theo:joint_density} Let $R^*$ and $H^*$ denote the distance and height, respectively, of the specific building that determines the global maximum blockage elevation angle over all directions. The joint PDF of $(R^*, H^*)$ is given by
\begin{align}\label{eq:theojoint}
    &\mathbb{P}[R^* \in [r, r+\mathrm{d}r], H^* \in [h, h+\mathrm{d}h]] = 2\pi\lambda r f_H(h) \nonumber \\
    &\quad \times \exp\biggl(-2\pi\lambda \int_{0}^{\infty} \biggl[ 1 - F_H\biggl( \frac{h y}{r} \biggr) \biggr] y \,\mathrm{d}y\biggr) \,\mathrm{d}r \,\mathrm{d}h.
\end{align}

\begin{IEEEproof}
A building located at a distance $r$ and with height $h$ creates the global maximum elevation angle when the following two conditions are satisfied:
\begin{enumerate}
\item A building exists in the infinitesimal area $r \mathrm{d}r \mathrm{d}\theta$ with height in $[h, h+\mathrm{d}h]$. The probability of this occurrence is $\lambda (2\pi r \mathrm{d}r) f_H(h) \mathrm{d}h$.
\item No other building in the entire spatial domain $\mathbb{R}^2$ generates an elevation angle larger than $\phi = \arctan(h/r)$. 
\end{enumerate}
The second condition corresponds to the probability that the number of points falling into the blocking region is zero, where the blocking region $\mathcal{S}_{\text{block}}$ is the set of points $(r', \theta', h')$ where $\frac{h'}{r'} > \frac{h}{r}$. 

The mean number of points in this region is calculated as
\begin{align}
    \Lambda(\mathcal{S}_{\text{block}}) 
    &= \lambda \int_0^{2\pi} \int_0^\infty \mathbb{P}\biggl[ H' > r' \frac{h}{r} \biggr] r' \,\mathrm{d}r' \,\mathrm{d}\theta' \nonumber \\
    &= 2\pi\lambda \int_0^\infty \biggl[ 1 - F_H\biggl( \frac{h r'}{r} \biggr) \biggr] r' \,\mathrm{d}r'.\nonumber
\end{align}
The probability that no points exist in this region is $\exp(-\Lambda(\mathcal{S}_{\text{block}}))$. Combining this with the density of $(r,h)$ yields \eqref{eq:theojoint}. Note that, in the integral, we replaced $r'$ with $y$ to match the theorem statement.
\end{IEEEproof}

\example[Exponential height distribution]
Assume that the building height distribution is exponential with rate parameter $\mu$, i.e., $F_H(h) = 1 - e^{-\mu h}$.
Then, from Theorem~\ref{theo:joint_density}, the joint density simplifies to
\begin{align}
    2\pi\lambda \mu r \exp\biggl(-\frac{2\lambda\pi r^2}{h^2\mu^2} - \mu h\biggr) \,\mathrm{d}r \,\mathrm{d}h.\nonumber
\end{align}
From this joint distribution, we can derive the marginal distributions. Integrating \eqref{eq:theojoint} with respect to $r$ yields the marginal PDF of the height $H$ as
\begin{align}
    g(h) &= \frac{\mu^3}{2} h^2 e^{-\mu h}, \quad h \ge 0.\nonumber
\end{align}
This corresponds to a Gamma distribution with shape parameter $k=3$ and rate parameter $\theta=\mu$, i.e., $H \sim \text{Gamma}(3, \mu)$. 

Similarly, the marginal PDF of the distance $R$ becomes
\begin{align}
    k(r) &= 2\lambda \sqrt{\pi} \, r \, G^{0,3}_{3,0}
    \left(
    \begin{array}{c}
    - \\
    0, \tfrac12, 1
    \end{array}
    \middle\vert
    \frac{\lambda \pi r^2}{2}
    \right),\nonumber
\end{align}
where $G^{\cdot,\cdot}_{\cdot,\cdot}(\cdot|\cdot)$ denotes the Meijer $G$-function \cite{beals2013meijer}. The corresponding means are
\begin{align}
    \mathbb{E}[H] = \frac{3}{\mu},~~    \mathbb{E}[R] &= \frac{3}{2\sqrt{2\lambda}}.\nonumber
\end{align}

We can observe that the average height of the building that creates the largest blockage elevation angle is $\mathbb{E}[H] = 3/\mu$. This indicates that high buildings are more likely to determine the global maximum blockage angle because they are visible and dominant over longer distances. Thus, the skyline is not determined by the average building but by outliers that are, on average, three times taller than the typical building.

\vspace{-1em}
\subsection{Joint Spatial Statistics for Angular Separation}\label{subsec:3-3}
\noindent\textbf{Joint Spatial Statistics for Angular Separation:}
In a dense urban environment, buildings can block the signals from LEO satellites to the ground user. To mitigate this, connecting to multiple satellites that are spatially separated is essential. However, the phenomenon of blocking of two adjacent satellites is not independent. If a large building blocks one satellite, nearby satellites are also likely to be blocked by the same obstacle. Therefore, to design a robust dual-connectivity strategy, it is essential to quantify the spatial correlation of the blocking over the azimuth angle. For this, we analyze the probability of joint outage of two satellites separated by $\Delta \psi$ in azimuth angle. Using the joint statistics of the skyline process derived in Theorem \ref{theo:joint_CDF_omega_1_psi_1_psi_2} of Section \ref{sec:sec3}, the probability that two satellites will be blocked at the same time can be evaluated. The probability of disconnection of the dual link $P_{\text{out}}^{\text{dual}}(\theta, \Delta \psi)$ is defined as the probability that the elevation angles of blockage of both satellites are lower (i.e., blocked) than the skyline.

We now examine the joint distribution of $\omega_1(\psi_1)$ and $\omega_1(\psi_2)$ for two directions $\psi_1$ and $\psi_2$ separated by an angular difference $\Delta_\psi$ where $\psi_1,\psi_2\in(-\pi,\pi]$. More specifically, let $\Delta_{\psi}=\min(|\psi_1-\psi_2|,2\pi-|\psi_1-\psi_2|)$.
Buildings whose centroids lie in the overlapping region $\mathcal{A}(\psi_1)\cap \mathcal{A}(\psi_2)$ can block both directions, $\psi_1$ and $\psi_2$, and this region is
\begin{equation}\footnotesize \label{eq:intersection}
\begin{aligned}
    &\mathcal{A}(\psi_1) \cap \mathcal{A}(\psi_2) = \\
    &\begin{cases} 
        \{ (r, \alpha) : \alpha \in [-\pi, \pi] \}, & \text{if } r < \frac{l}{2\pi} \\
        \{ (r, \alpha) : \alpha \in [\psi \pm \frac{l}{2r}] \bmod 2\pi \}, & \text{if } r \geq \frac{l}{2\pi}, \Delta_\psi \leq \frac{l}{r} \\
        \emptyset, & \text{otherwise}
    \end{cases}\nonumber
\end{aligned}
\end{equation}
Further, 
\begin{equation}\footnotesize \label{eq:difference}
\begin{aligned}
    &\mathcal{A}(\psi_1) \setminus \mathcal{A}(\psi_2) = \\
    &\begin{cases}
        \emptyset, & \text{if } r < \frac{l}{2\pi} \\
        \{ (r, \alpha) : \alpha \in (\psi_1 \pm \tfrac{l}{2r}) \setminus (\psi_2 \pm \tfrac{l}{2r}) \}, & \text{if } r \geq \frac{l}{2\pi}
    \end{cases}\nonumber
\end{aligned}
\end{equation}

\theorem\label{theo:joint_CDF_omega_1_psi_1_psi_2} The joint CDF of $\omega_1(\psi_1)$ and $\omega_1(\psi_2)$ is 
\begin{align}\label{eq:joint_cdf_psi1_psi2}
    &\mathbb{P}[\omega_1(\psi_1)<\phi_1 \cap \omega_1(\psi_2)<\phi_2] \nonumber \\
    &= e^{-\lambda\int_{\mathcal{A}(\psi_1)\cap\mathcal{A}(\psi_2)} [1 - F_H(r \min(\tan\phi_1,\tan\phi_2))] d\mathbf{x}} \nonumber \\
    &\quad \times e^{-\lambda\int_{\mathcal{A}(\psi_1)\setminus\mathcal{A}(\psi_2)} [1 - F_H(r \tan\phi_1)] d\mathbf{x}} \nonumber \\
    &\quad \times e^{-\lambda\int_{\mathcal{A}(\psi_2)\setminus\mathcal{A}(\psi_1)} [1 - F_H(r \tan\phi_2)] d\mathbf{x}},
\end{align}
where $\mathbf{x} = (r, \alpha)$ denotes the polar coordinates of a building centroid, and $r$ is the distance from the user.
\begin{IEEEproof}
The joint CDF of the tangents of $\omega_1(\psi_1)$ and $\omega_1(\psi_2)$ is
\begin{align}\label{eq:proof_eq:joint_cdf_psi1_psi2}
    &\mathbb{P}[\tan\omega_1(\psi_1)<\tau_1,\tan\omega_1(\psi_2)<\tau_2] \nonumber \\
    &= \mathbb{P}\biggl[\sup_{\mathbf{x}_i\in\Phi\cap(\mathcal{A}_1\cap\mathcal{A}_2)}\frac{h_i}{r_i}\leq \min(\tau_1,\tau_2)\biggr] \nonumber \\
    & \times \mathbb{P}\biggl[\sup_{\mathbf{x}_i\in\Phi\cap(\mathcal{A}_1\setminus\mathcal{A}_2)}\frac{h_i}{r_i}\leq \tau_1\biggr] \mathbb{P}\biggl[\sup_{\mathbf{x}_i\in\Phi\cap(\mathcal{A}_2\setminus\mathcal{A}_1)}\frac{h_i}{r_i}\leq \tau_2\biggr].
\end{align}
As in Theorem \ref{theo:CDF_omega1_psi}, we can obtain the joint CDF of $\omega_1(\psi_1)$ and $\omega_1(\psi_2)$ as \eqref{eq:joint_cdf_psi1_psi2} by applying the integration region of $\mathcal{A}(\psi_1)\cap\mathcal{A}(\psi_2)$, $\mathcal{A}(\psi_1)\setminus\mathcal{A}(\psi_2)$, $\mathcal{A}(\psi_2)\setminus\mathcal{A}(\psi_1)$ in \eqref{eq:proof_eq:joint_cdf_psi1_psi2}, respectively.
\end{IEEEproof}

Equation~\eqref{eq:joint_cdf_psi1_psi2} implies that the joint probability of the blockage elevation angles in two directions decreases as the angular separation $\Delta_\psi$ increases. The spatial correlation between $\omega_1(\psi_1)$ and $\omega_1(\psi_2)$ is decreased, because fewer buildings simultaneously affect both viewing angles. In the extreme case where $\Delta_\psi \to 0$, the overlap of $\mathcal{A}(\psi_1)$ and $\mathcal{A}(\psi_2)$ becomes complete, and the two angles are nearly identical.

\remark When $\psi_1,\psi_2=\psi$ and $\phi_1,\phi_2=\phi$, \eqref{eq:joint_cdf_psi1_psi2} reduces to \eqref{eq:CDF_omega_1}.

\example When the building heights follow an exponential distribution with mean $1/\mu$, the joint CDF is simplified as
\begin{align}\label{eq:joint_cdf_exp_simple}
    &\mathbb{P}[\omega_1(\psi_1)<\phi, \omega_1(\psi_2)<\phi]= \exp \Biggl( - \lambda \biggl[ {\int_{\mathcal{S}_{\cap}} e^{-\mu r \tan\phi} d\mathbf{x}} \nonumber \\
    &\quad + {\int_{\mathcal{S}_{1\setminus 2}} e^{-\mu r \tan\phi} d\mathbf{x}} + {\int_{\mathcal{S}_{2\setminus 1}} e^{-\mu r \tan\phi} d\mathbf{x}} \biggr] \Biggr),
\end{align}
where $\mathcal{S}_{\cap} = \mathcal{A}(\psi_1)\cap\mathcal{A}(\psi_2)$, and $\mathcal{S}_{i\setminus j} = \mathcal{A}(\psi_i)\setminus\mathcal{A}(\psi_j)$.

Although a closed-form expression can be derived by evaluating these integrals over the polygonal intersection areas, it results in a lengthy formula. Instead, we numerically evaluate \eqref{eq:joint_cdf_exp_simple} to analyze the correlation behavior.

\begin{figure}[t]
\begin{center}
\includegraphics[width=0.78\linewidth]{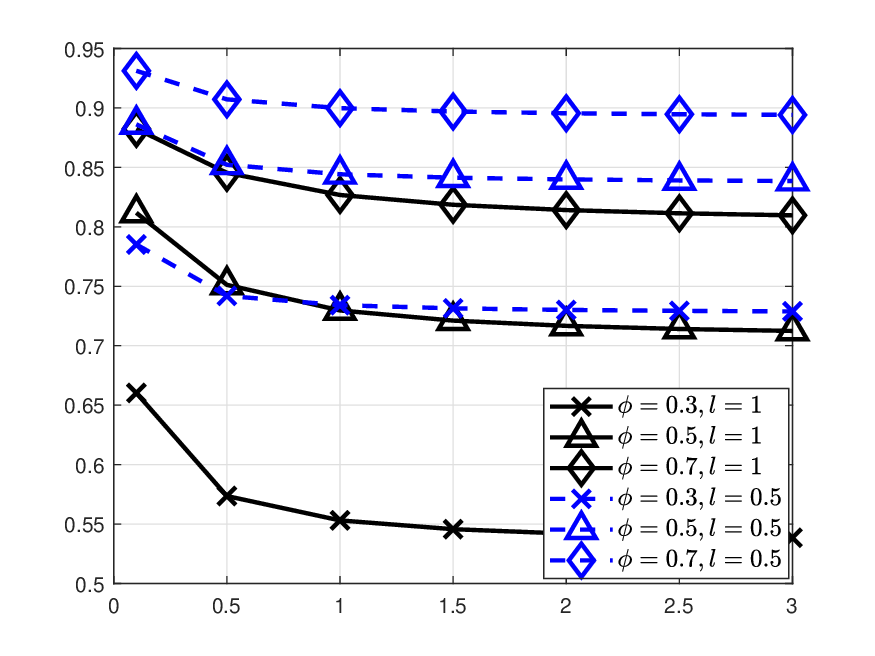}
\end{center}
\vspace{-0.5em} 
\caption{Joint CDF when varing $\Delta_{\psi}$ when $\mu = 1$, and $\lambda = 0.1$.} 
\label{fig:joint_CDF_psi}
\vspace{-1em}\end{figure}

Figure~\ref{fig:joint_CDF_psi} illustrates the joint CDF of $\omega_1(\cdot)$ for different $\Delta_\psi$, when the building heights follow an exponential distribution with $\mu = 1$ and when $\lambda = 0.1$. The joint CDF decreases as $\Delta_\psi$ increases, confirming the spatial decorrelation.

Regarding the building width $l$, while longer buildings intuitively increase the spatial correlation between directions, the figure shows that the joint CDF decreases since the increase in the individual blockage probability dominates the correlation effect. Even though blockage events are more correlated, the likelihood of simultaneously observing low blockage angles in both directions diminishes.

\noindent\textbf{Probability of $\boldsymbol{\omega_1(\psi_1) = \omega_1(\psi_2)}$:}
In the following theorem, we derive the probability that the blockage elevation angles in two different directions, $\omega_1(\psi_1)$ and $\omega_1(\psi_2)$ are exactly the same. Such an event occurs when the same building simultaneously creates $\omega_1(\cdot)$ in both directions. 

\theorem\label{theo:prob_omega_1_psi_1_psi_2_equal} Given $\psi_1$ and $\psi_2$, the probability that $\omega_1(\psi_1)$ and $\omega_1(\psi_2)$ are equal is 
{\begin{align}\label{eq:equalpsi1psi2}
&\int_{0}^{\frac{\pi}{2}} \Biggl[ \left( \frac{\mathrm{d}}{\mathrm{d}\phi} \exp\left( -\lambda \int_{\mathcal{A}(\psi_1)\cap\mathcal{A}(\psi_2)} 1-F_H(\|\mathbf{x}\|\tan\phi) \mathrm{d}\mathbf{x} \right) \right) \nonumber \\
&\quad \times \exp\left( -\lambda \int_{\mathcal{A}(\psi_1)\triangle\mathcal{A}(\psi_2)} 1-F_H(\|\mathbf{x}\|\tan\phi) \mathrm{d}\mathbf{x} \right) \Biggr] \mathrm{d}\phi.
\end{align}}
\begin{IEEEproof}
    Since the locations of buildings are Poisson, almost surely $\omega_1(\psi_1) = \omega_1(\psi_2)$ occurs if and only if the blockage elevation angles in both directions are determined by the same building. This condition requires two simultaneous events for a given angle $\phi$, where: 
\begin{enumerate}
\item The dominant building which makes the blockage angle $\phi$ must be located in the overlapping region $\mathcal{A}(\psi_1) \cap \mathcal{A}(\psi_2)$. The probability density for this is represented by the derivative term in \eqref{eq:equalpsi1psi2}, where $\|\mathbf{x}\|$ denotes the Euclidean distance from the user to point $\mathbf{x}$.
\item Other buildings in $\mathcal{A}(\psi_1) \triangle \mathcal{A}(\psi_2)$ do not create a blockage elevation angle larger than $\phi$. 
\end{enumerate}
    Integrating the product of these terms over all possible angles $\phi \in [0, \pi/2)$ yields the total probability.
\end{IEEEproof}

This theorem quantifies how the geometric overlap affects the directional persistence of the blockage. When $\Delta_\psi$ is small, the two regions share a larger spatial region, the probability that the same building obstructs both directions increases. 

\subsection{Second-Order Statistics and Spectral Analysis} 
Based on the spatial correlation derived in Section~\ref{subsec:3-3}, we now extend our analysis to the frequency domain to characterize the roughness or rate of change of $\omega_1(\psi)$. This spectral analysis provides insights into how quickly blockage conditions change as a user rotates.

\subsubsection{Auto-Correlation Function of $\omega_1(\psi)$}\label{subsec:3-4-1}
\mbox{}\\
\noindent\textbf{Second-Order Statistics (ACF):}
$\omega_1(\psi)$ is a rotation invariant stochastic process with respect to the azimuth angle $\psi$. To quantify the correlation of $\omega_1(\psi)$ over the entire angular region, we define the auto-correlation function (ACF) of $\omega_1(\psi)$, denoted as $R_{\omega_1}(\delta)$.

In Theorem~\ref{theo:joint_CDF_omega_1_psi_1_psi_2}, we derived the joint CDF which represents the probability that the blockage elevation angles at two directions are lower than $\phi_1$ and $\phi_2$, respectively. Let us denote this joint CDF when two directions are separated by an angular distance $\delta$ as $F_{\omega_1}(\phi_1, \phi_2; \delta)$. Since the azimuth angle $\psi$ is defined on a circle, the angular positions are considered modulo $2\pi$. In other words,
\begin{align}
F_{\omega_1}(\phi_1, \phi_2; \delta) = \mathbb{P}[\omega_1(\psi) \leq \phi_1, \omega_1((\psi+\delta) \mod{2\pi}) \leq \phi_2].\nonumber
\end{align}
Then, the joint PDF, $f_{\omega_1}(\phi_1, \phi_2; \delta)$ becomes
\begin{align}
    f_{\omega_1}(\phi_1, \phi_2; \delta) = \frac{\partial^2}{\partial \phi_1 \,\partial \phi_2} \, F_{\omega_1}(\phi_1, \phi_2; \delta).\nonumber
\end{align}
Taking the circular stationarity into account, the ACF for an angular separation $\delta$ is obtained as
\begin{align}
    R_{\omega_1}(\delta) 
    &= \mathbb{E}\bigl[ \omega_1(\psi) \, \omega_1((\psi+\delta) \bmod 2\pi) \bigr] \nonumber \\
    &= \int_{0}^{\pi/2} \int_{0}^{\pi/2} \phi_1 \phi_2 f_{\omega_1}(\phi_1, \phi_2; \delta) \,\mathrm{d}\phi_1 \,\mathrm{d}\phi_2.\nonumber
\end{align}

\begin{figure}[!htbp]
    \centering
    \begin{subfigure}{1.0\columnwidth}
        \centering
        \includegraphics[width=0.78\linewidth]{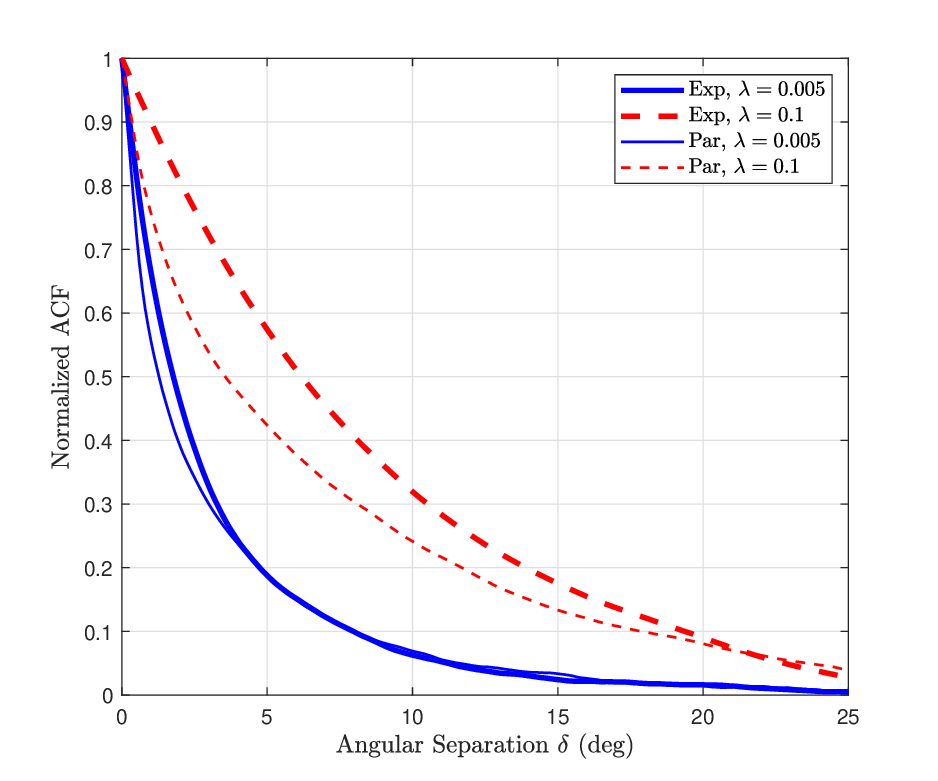}
        \caption{ACF of $\omega_1(\psi)$.}
        \label{fig:ACF}
    \end{subfigure}

    \begin{subfigure}{1.0\columnwidth}
        \centering
        \includegraphics[width=0.78\linewidth]{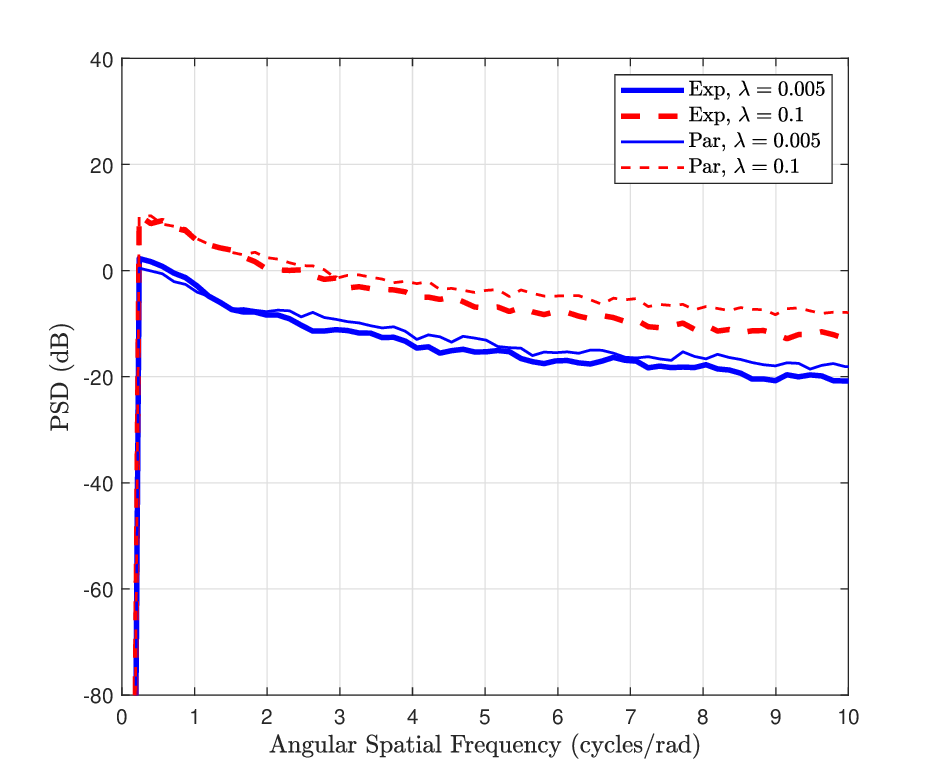}
        \caption{PSD of $\omega_1(\psi)$.}
        \label{fig:PSD}
    \end{subfigure}
    \caption{ACF and PSD of $\omega_1(\psi)$.}
    \label{fig:combined_acf_psd}
\vspace{-1em}\end{figure}

Figure~\ref{fig:ACF} illustrates the ACF, $R_{\omega}(\delta)$ when $l=1$. Further, we set $\mu=1$, $\kappa=1.5$, and $s=\frac{1}{3}$ to make the average building height $1$. A slowly decaying ACF indicates a highly correlated environment over the angular separation. In dense environments ($\lambda=0.1$), the skyline remains correlated over larger angles, as buildings closer to the user occupy a wider field of view. Notably, the Pareto distribution exhibits a longer coherence angle than the Exponential distribution, as heavy-tailed heights create dominant, persistent blockages.

\subsubsection{Power Spectral Density of $\omega_1(\psi)$}\label{subsec:3-4-2}
\mbox{}\\
\noindent\textbf{Spectral Analysis (PSD):}
We apply the Wiener-Khinchin theorem to this ACF to analyze the fluctuation characteristics of $\omega_1(\psi)$ in the frequency domain. The power spectral density (PSD) of $\omega_1(\psi)$, $S_{\omega_1}(f)$, is obtained by applying the Fourier transform to the ACF. So,
\begin{align}
S_{\omega_1}(f) = \int_{-\infty}^{\infty} R_{\omega_1}(\xi) e^{-j 2\pi f \delta} \mathrm{d}\xi,\nonumber
\end{align}
where $f$ denotes the angular spatial frequency.

Figure~\ref{fig:PSD} shows the PSD, $S_{\omega}(f)$ when $l=1$. Further, we set $\mu=1$, $\kappa=1.5$, and $s=\frac{1}{3}$ to make the average building height $1$ as in the ACF part. In dense environments ($\lambda=0.1$, we observe higher PSD, indicating larger variance in skyline fluctuations. Especially, the Pareto distribution shows higher energy in the low-frequency region than in the Exponential case, indicating that the skyline is characterized by larger, more continuous blockage due to extremely tall buildings.


\section{Application to LEO Satellite Network Design}\label{sec:application}
This section analyzes the reliability of LEO satellite networks in urban environments using the statistical characteristics of the skyline process derived in Section~\ref{sec:sec3}. Specifically, we investigate the impact of the minimum elevation mask angle on network availability, quantified by the outage probability and satellite diversity gain.

\subsection{Satellite Diversity and Outage Probability}\label{subsec:application1}
\subsubsection{System Model for Section~\ref{subsec:application1}}

We consider LEO satellites to be modeled as a homogeneous PPP with intensity $\lambda_{sat}$ [satellites/km$^2$] on the sphere of radius $R_{orb} = R_E + h_{sat}$, where $R_E$ is the Earth's radius and $h_{sat}$ is the satellite altitude. Let the total number of satellites on the sphere with radius $R_{orb}$ be $N_{total}$. Then, $N_{total}$ follows a Poisson distribution with mean $\mathbb{E}[N_{total}]=\bar{N}_{total}=4\pi R_{orb}^2 \lambda_{sat}$.

For a satellite, from the view of a ground user, the azimuth angle $\psi$ is uniformly distributed in $[0, 2\pi)$, and the PDF of the elevation angle $\theta$ for a potentially visible satellite is given by
\begin{align}\label{eq:exact_PDF}
    f_{\Theta}(\theta) = \frac{\cos \theta}{(1-k)\sqrt{1 - k^2 \cos^2 \theta}}  \biggl( \sqrt{1 - k^2 \cos^2 \theta} - k \sin \theta \biggr)^2, 
\end{align}
where $k = R_E / (R_E + h_{sat})$, for $\theta \in [0, \pi/2]$. The detailed proof is provided in Appendix~\ref{appen:elevation_pdf}. 

A satellite is visible to the user if 1) its elevation angle from the user, $\theta$, exceeds the minimum elevation mask angle $\theta_{\min}$, and 2) the path between the user and the satellite is not obstructed by buildings. The latter is satisfied when $\theta > \omega_1(\psi)$, where $\omega_1(\cdot)$ is defined in Section~\ref{sec:sec2}. Consequently, a satellite is visible if $\theta > \max(\theta_{\min}, \omega_1(\psi))$. 

In our numerical evaluation, we assume that the mean total number of satellites on the sphere with radius $R_{orb}$ is $\bar{N}_{total} = 10,000$, where $h_{sat} = 500$km. Buildings are assumed to have a width of $l=50$ m, with heights following an exponential distribution.

\subsubsection{Mean Number of Visible Satellites ($\bar{N}_{vis}(\theta_{\min})$)}
The mean number of visible satellites given $\theta_{\min}$, $\bar{N}_{vis}(\theta_{\min})$, is calculated by multiplying $\bar{N}_{total}$ by the probability that a single satellite at an arbitrary location is observable by the user, $P_{vis}(\theta_{\min})$. The latter is obtained by integrating the LoS probability over the valid sky region and is given by
\begin{align}
    P_{vis}(\theta_{\min}) = \int_{\theta_{\min}}^{\pi/2} F_{exp}(\theta){f_{\Theta}(\theta)} \mathrm{d}\theta,\nonumber
\end{align}
where $F_{exp}(\theta)$ and $f_{\Theta}(\theta)$ are the CDF $\omega_1(\psi)$ under the exponential distribution with parameter $\mu$ in \eqref{eq:ex_CDF_omega1_exp} and the PDF of the satellite elevation angle in \eqref{eq:exact_PDF}, respectively. Note that the CDF of $\omega_1(\psi)$ is given in Theorem~\ref{theo:CDF_omega1_psi}. So, the mean number of visible satellites is given by
\begin{align}
    \bar{N}_{vis}(\theta_{\min}) = \bar{N}_{total} P_{vis}(\theta_{\min}).\nonumber
\end{align}

\begin{figure}[!htbp]
\vspace{-1em}
    \centering
    \begin{subfigure}{1.0\columnwidth}
        \centering
        \includegraphics[width=0.78\linewidth]{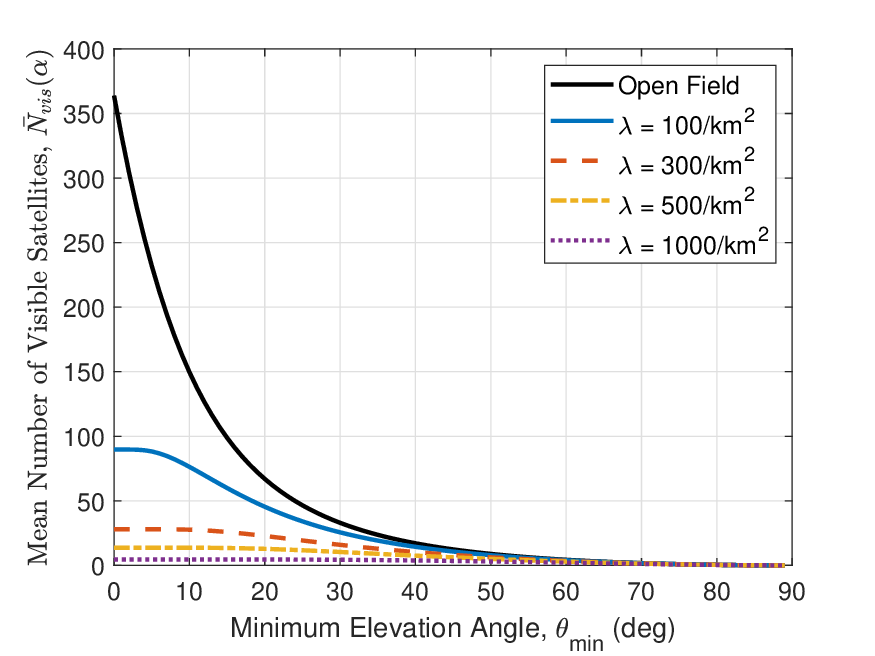}
        \caption{Impact of building density $\lambda$ with $\mu^{-1}=50$m.}
        \label{fig:visibility_lambda}
    \end{subfigure}

    \begin{subfigure}{1.0\columnwidth}
        \centering
        \includegraphics[width=0.78\linewidth]{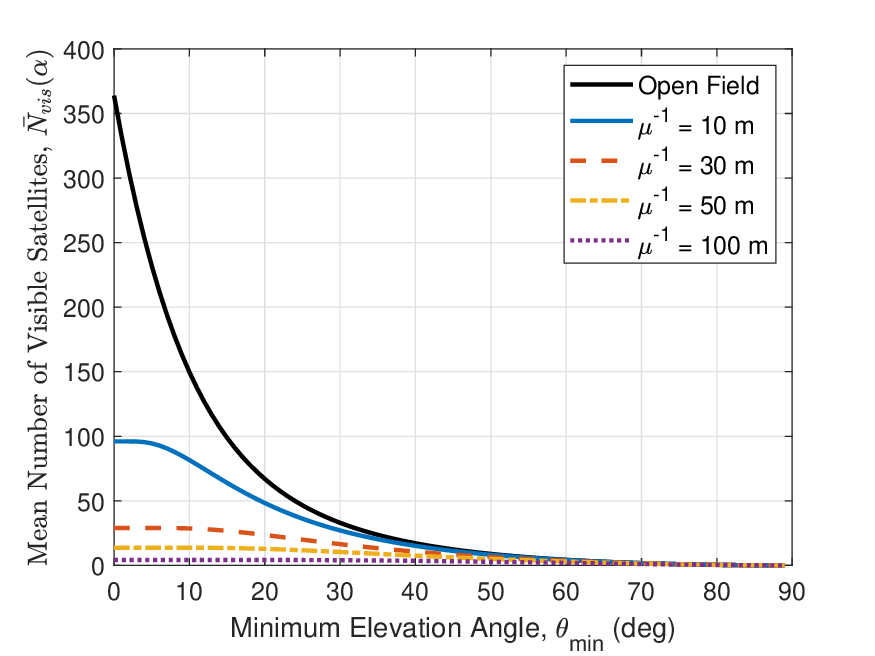}
        \caption{Impact of building height $\mu^{-1}$ with $\lambda=500/\text{km}^2$.}
        \label{fig:visibility_mu}
    \end{subfigure}

    \caption{Effects of building parameters on the average number of visible satellites $\bar{N}_{vis}(\theta_{\min})$.}
    \label{fig:visibility_comparison}
\vspace{-1em}\end{figure}

Figure~\ref{fig:visibility_comparison} illustrates $\bar{N}_{vis}(\theta_{\min})$ with different $\theta_{\min}$, $\lambda$, and $\mu^{-1}$. In an open field case ($\lambda=0$ or $\mu^{-1}=0$), when $\theta_{\min}=0$, approximately $360$ satellites are observed by the user, but $98\%$ of these satellites are blocked by buildings as $\lambda$ and $\mu^{-1}$ increases to $1000/\text{km}^2$ and $100$m, respectively. The curves in the two figures appear to scale similarly with the ratio $\lambda/\mu$, but, in fact, the non-linear dependency on $\mu$ within the exponential term of \eqref{eq:ex_CDF_omega1_exp} makes these curves not identical. A smaller $\theta_{\min}$ makes the observable sky region wider, but degrades signal-to-noise ratios (SNRs). So, choosing an appropriate $\theta_{\min}$ to balance between the link reliability and the total number of visible satellites is needed.

\subsubsection{Outage Probability ($P_{out}(\theta_{\min})$)}

We define the outage probability as the probability that the user cannot observe any satellites given a minimum elevation mask angle $\theta_{\min}$. 
If we assume that blockage events for different satellites are independent (the independent thinning approximation), the number of visible satellites follows a Poisson distribution with mean $\bar{N}_{vis}(\theta_{\min})$. In this case, the outage probability is:
\begin{align}\label{eq:Pout_ind}
&P_{out}^{ind}(\theta_{\min}) = \exp\left( - \bar{N}_{vis}(\theta_{\min}) \right) \nonumber\\
&= \exp\left( - \bar{N}_{total} \int_{\theta_{\min}}^{\pi/2} \mathbb{P}[\omega_1(\psi) < \theta] f_{\Theta}(\theta) \mathrm{d}\theta \right).
\end{align}
However, this approximation does not consider the positive spatial correlation. To investigate this correlation, we define $I(\gamma, r, h, \psi)$ for a building with azimuth $\gamma$, height $h$ and distance $r$ as
\begin{align}
I(\gamma, r, h, \psi) = 
\begin{cases} 
\theta, & \text{if } \psi \in [\gamma - \frac{\nu}{2}, \gamma + \frac{\nu}{2}] \pmod{2\pi} \\
0, & \text{otherwise, }
\end{cases}\nonumber
\end{align}
where the blockage elevation angle is $\theta = \arctan(h/r)$ and the angular half-width is $\nu =\min( \frac{\ell}{2r}, \pi)$. 

For all satellites in \textcolor{black}{$\Phi_{sat}$}, the blockage event $A_i$ at azimuth $\psi_i$ is defined as
\begin{align}
A_i = \left\{ \sup_{j \in \Phi_b} I(\gamma_j, r_j, h_j, \psi_i) > \theta_{\min} \right\},\nonumber
\end{align}
where $\Phi_b$ is the building process.

To analyze the correlation between any two directions $\psi_i$ and $\psi_{i'}$, which are the azimuth angles of the $i$-th and $i'$-th satellites, we sort the buildings into four categories as
\begin{itemize}
    \item $\mathcal{D}_i$: Buildings which cross $\psi_i$ but not $\psi_{i'}$.
    \item $\mathcal{D}_{i'}$: Buildings which cross $\psi_{i'}$ but not $\psi_i$.
    \item $\mathcal{B}$: Buildings which cross both $\psi_i$ and $\psi_{i'}$.
    \item $\mathcal{O}$: Buildings that block neither.
\end{itemize}

The joint blockage probability $\mathbb{P}(A_i \cap A_{i'})$ is determined by these independent sets. Let $E(\mathcal{X})$ be the event that at least one building in the set $\mathcal{X}$ blocks the satellite, and then $A_i = E(\mathcal{D}_i) \cup E(\mathcal{B})$ and $A_{i'} = E(\mathcal{D}_{i'}) \cup E(\mathcal{B})$.
\begin{align}
    &\mathbb{P}(A_i \cap A_{i'} \mid \textcolor{black}{\Phi_{sat}}) = \mathbb{P}\bigl( (E(\mathcal{D}_i) \cup E(\mathcal{B})) \cap (E(\mathcal{D}_{i'}) \cup E(\mathcal{B})) \bigr) \nonumber \\
    &= \mathbb{P}(E(\mathcal{D}_i) \cap E(\mathcal{D}_{i'})) + \mathbb{P}(E(\mathcal{B})) \ge \mathbb{P}(A_i)\mathbb{P}(A_{i'}).\nonumber
\end{align}
By applying the FKG inequality over $\mathcal{B}$ \cite{holley1974remarks}, we obtain the relationship of outage probabilities under the independent approximation and the real environment as
\begin{align}
    P_{out}(\theta_{\min}) 
    &= \textcolor{black}{\mathbb{E}_{\Phi_{sat}}} \Biggl[ \mathbb{P} \Biggl( \textcolor{black}{\bigcap_{x_i \in \Phi_{sat}}} A_i \textcolor{black}{\Bigm| \Phi_{sat}} \Biggr) \Biggr] \nonumber \\
    &\ge \textcolor{black}{\mathbb{E}_{\Phi_{sat}}} \Biggl[ \textcolor{black}{\prod_{x_i \in \Phi_{sat}}} \mathbb{P}(A_i) \Biggr] = P_{out}^{ind}(\theta_{\min}).\nonumber
\end{align}
So, the approximation in \eqref{eq:Pout_ind} in fact provides a lower bound on the outage probability.

\begin{figure}[!htbp]
\vspace{-1em}
    \centering
    \begin{subfigure}{1.0\columnwidth}
        \centering
        \includegraphics[width=0.78\linewidth]{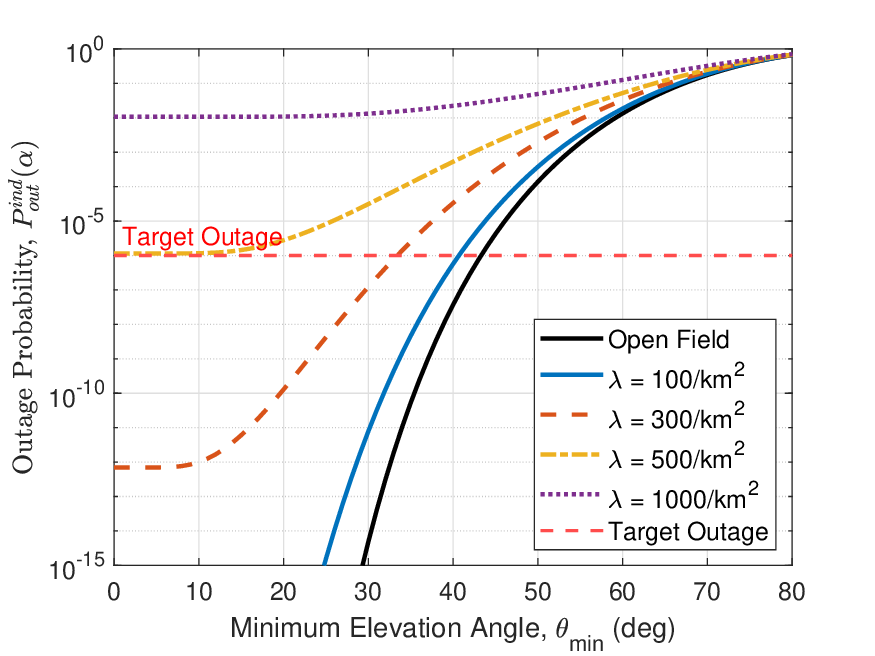}
        \caption{Varying building density $\lambda$ with $\mu^{-1}=50$m.}
        \label{fig:outage_lambda}
    \end{subfigure}       
    \begin{subfigure}{1.0\columnwidth}
        \centering
        \includegraphics[width=0.78\linewidth]{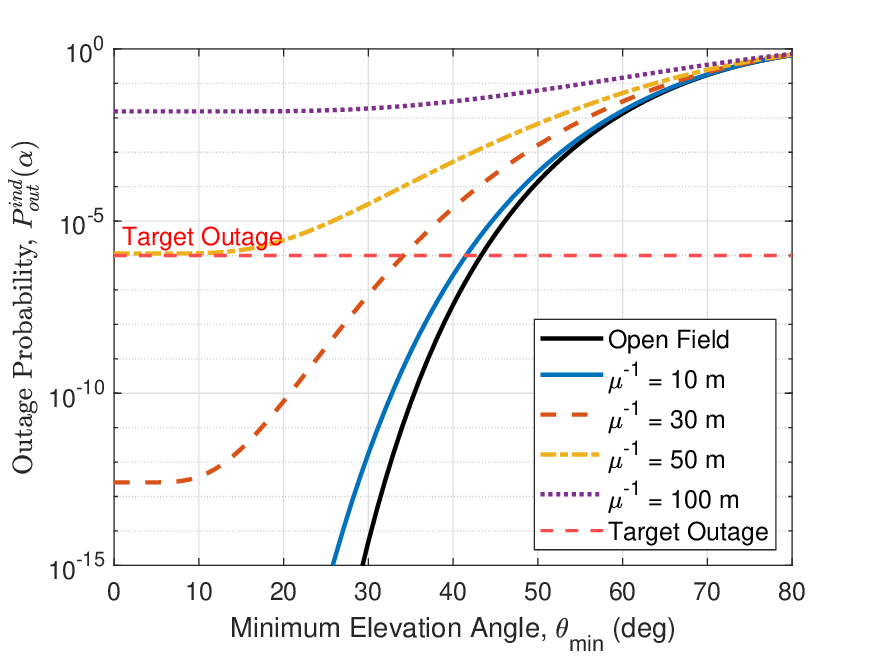}
        \caption{Varying building height $\mu^{-1}$ with $\lambda=500/\text{km}^2$.}
        \label{fig:outage_mu}
    \end{subfigure}
    \caption{Impact of building parameters on the outage probability $P_{out}(\theta_{\min})$. The dashed horizontal line represents the target outage probability of $10^{-6}$.}
    \label{fig:outage_comparison}
\vspace{-1em}\end{figure}

Figure~\ref{fig:outage_comparison} depicts the relationship between $\lambda$, $\theta_{\min}$, and $P^{ind}_{out}(\theta_{\min})$, with a target outage of $10^{-6}$. As $\lambda$ increases or $\mu$ decreases, the outage probability curve shifts upward. When $\lambda$ is larger than $500km/m^2$ under $\mu^{-1}=50m$ or $\mu^{-1}$ is larger than $50m$ under $\lambda=500/km^2$, decreasing $\theta_{\min}$ cannot meet the target outage $10^{-6}$. As in Figure~\ref{fig:visibility_comparison}, the two figures in Figure~\ref{fig:outage_comparison} appear to scale similarly with the ratio $\lambda/\mu$, but these are not identical due to the exponential term in \eqref{eq:ex_CDF_omega1_exp}. While a higher value of $\theta_{\min}$ is preferred to reduce atmospheric attenuation and terrestrial interference, selecting the appropriate $\theta_{\min}$ is crucial for reliable connectivity.

\subsection{Impact of Angular Separation on Diversity Gain}

In a dense urban environment, buildings frequently block signals from LEO satellites to ground users. To mitigate this, connecting to multiple spatially separated satellites is essential. However, the blockage events of two adjacent satellites are not independent. Therefore, it is crucial to quantify the spatial correlation of blockage relative to the azimuth angle.

We analyze the joint outage probability of two satellites located at the same elevation angle $\theta$, but separated by an azimuth angle $\Delta \psi$, and denote the probability that the blockage elevation angles in directions separated by $\Delta \psi$ are larger than $\theta$ by $P_{\text{out}}^{\text{dual}}(\theta, \Delta \psi)$ which is given by 
\begin{align}
    &P_{\text{out}}^{\text{dual}}(\theta, \Delta \psi) 
    = \mathbb{P}\bigl[ \omega_1(\psi) > \theta \cap \omega_1(\psi+\Delta \psi) > \theta \bigr] \nonumber \\
    &= 1 - 2\mathbb{P}\bigl[\omega_1(\psi) \le \theta\bigr]  + \mathbb{P}\bigl[\omega_1(\psi) \le \theta \cap \omega_1(\psi+\Delta \psi) \le \theta\bigr],\nonumber
\end{align}
since $\mathbb{P}[ \omega_1(\psi) \le \theta] = \mathbb{P}[ \omega_1(\psi+\Delta \psi) \le \theta]$ due to the stationarity property.

\begin{figure}[t]
    \centering
    \includegraphics[width=0.78\linewidth]{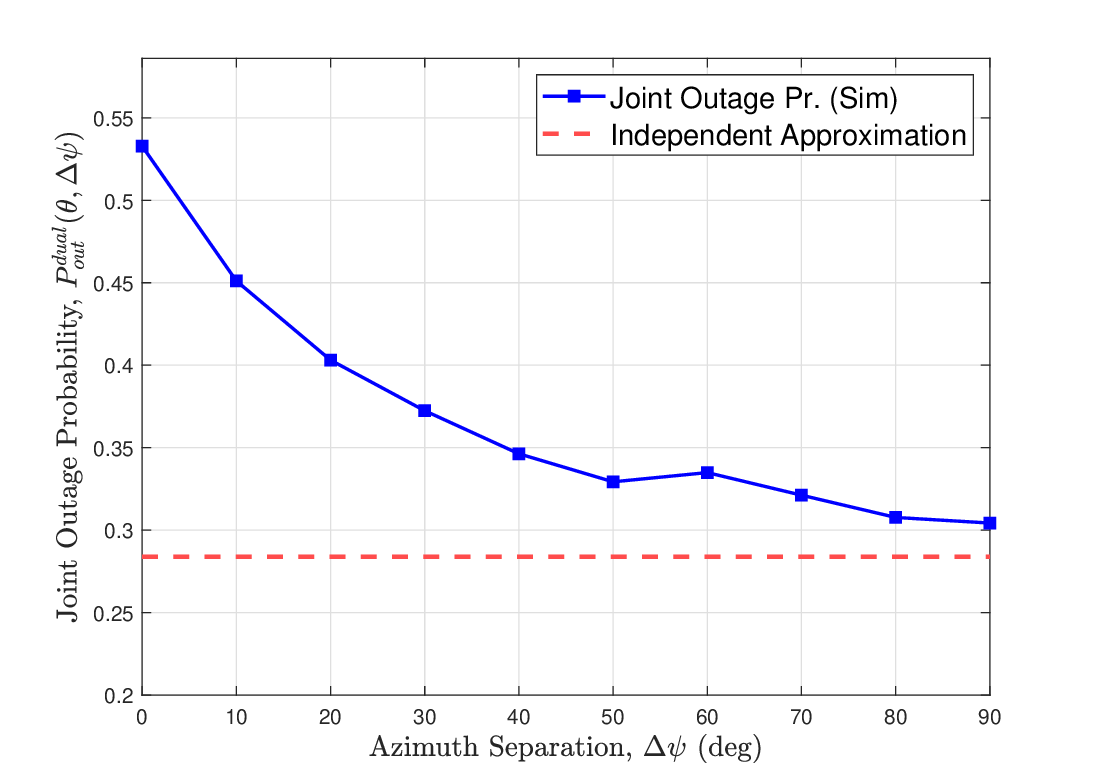}
    \caption{Joint outage probability with respect to $\Delta \psi$.}
    \label{fig:diversity}
\vspace{-1em}\end{figure}

Figure~\ref{fig:diversity} illustrates the effect of $\Delta \psi$ on $P_{\text{out}}^{\text{dual}}$. We compare the simulation results with the independent approximation, given by $(1 - \mathbb{P}[\omega_1(\psi) \le 45^{\circ}])^2$, where the blockage events are uncorrelated. For this evaluation, we assume that building heights follow an exponential distribution with $\mu^{-1}=30$\,m, width $l=25$\,m, and density $\lambda=1000/\text{km}^2$.

From these results, we derive three main observations:
\begin{itemize}
\item \textbf{Strong Correlation at Small Separation:} When $\Delta \psi$ is small (e.g., $<10^\circ$), the joint outage probability becomes close to the single-link outage probability, and this indicates a strong spatial correlation; if one satellite is blocked by a large building, its neighbors are highly likely to be blocked by the same obstacle.
\item \textbf{Correlation Resolution and Diversity Gain:} As $\Delta \psi$ increases, the joint outage probability decreases significantly. So, a spatial diversity gain can be expected once the angular spacing between satellites exceeds the average angular width of the obstructing buildings.
\item \textbf{Decorrelation Angle:} The curve converges to the independent curve at approximately $\Delta \psi \approx 40^\circ$, and this decorrelation angle provides critical design insights to maximize satellite diversity gain.
\end{itemize}

\section{Conclusion}\label{sec:conclusion} 

In this paper, we proposed a novel analytical framework for quantifying the sky visibility in 3D urban environments. By modeling the distribution of urban buildings using a 3D skyline process based on stochastic geometry, we derived closed-form formulas for the PDF and CDF of the blockage elevation angles. Furthermore, we extended the analysis beyond first-order statistics to investigate the spatial structure and spectral properties of the skyline process.

Our analysis provided both theoretical and practical insights to understand urban environments for communication between a ground user and NTN elements by characterizing the skyline process. We analyzed the distribution of the Skyline process by providing the closed-form expressions of blockage elevation angles, dominant obstacles, and joint spatial statistics for angular separation. We also provided second-order statistics of the skyline process under exponential and Pareto distributions. We observed that the statistics of the skyline process are highly dependent on urban environment parameters. 

Based on these findings, we provided practical insights to design the use of LEO satellite networks in urban environments by evaluating satellite diversity and outage probability. Based on these metrics, we discussed how to design the minimum elevation mask angle to compromise link reliability and link quality. We also studied the impact of angular separation on diversity gain. Due to the correlation between blocks, we need to treat the decorrelation angle as a critical design parameter for robust dual-connectivity.


\appendices
\section{Proof of the Ordering Relation between \eqref{eq:ex_CDF_omega1_exp} and \eqref{eq:ex_CDF_omega1_par}}\label{appen:eq1213}
\proposition $F^{\mbox{exp}}_{\omega_1}(\phi)$ and $F^{\mbox{par}}_{\omega_1}(\phi)$ are the CDFs of the elevation angle for Exponential and Pareto building height distributions, respectively. Under the condition that the mean heights are equal and $1<\kappa<2$, the following inequality holds for all elevation angles $\phi\in(0,\frac{\pi}{2})$, where
\begin{align}\label{eq:appen1}
    F^{\mbox{exp}}_{\omega_1}(\phi)>F^{\mbox{par}}_{\omega_1}(\phi).
\end{align}
\begin{IEEEproof}
    Proving \eqref{eq:appen1} is equivalent to proving 
    \begin{align}
        R(t) = \frac{\ln (F^{\mbox{exp}}_{\omega_1}(\phi))}{\ln (F^{\mbox{par}}_{\omega_1}(\phi))} = \frac{(1-e^{-t})t^{-2+\kappa}(\kappa-1)(2-\kappa)}{\left(\frac{\kappa-1}{\kappa}\right)^{\kappa}} <1,\nonumber
    \end{align}
    for all $t>0$. Here, we substituted $\tan\phi$ for $x$ and $\frac{lx}{2\pi}$ for $t$.
    
    \textbf{Step 1: }Asymptotic behavior ($t \rightarrow 0$ and $t\rightarrow \infty$)
    
    When $t\rightarrow 0$, using the approximation $t\simeq 1-e^{-t}$, the ratio scales as $R(t)\propto t^{\kappa-1}$. Since $\kappa>1$, $\lim_{t\rightarrow 0}R(t) = 0$. Also, when $t\rightarrow\infty$, $1-e^{-t}$ goes to 1. So, the ratio scales as $R(t)\propto t^{\kappa-2}$. Since $\kappa<2$, $\lim_{t\rightarrow \infty}R(t) = 0$.

    \textbf{Step 2: }Global Bound Analysis

    Since $R(t)$ is continuous for $t>0$, we analyze the first derivative of $\ln R(t)$ to find critical points, $t^*$<. 
    \begin{align}
        \frac{\mathrm{d}}{\mathrm{d}t}\ln R(t)=\frac{1}{e^t-1}-\frac{2-\kappa}{t} ~~\rightarrow~~ \frac{t^*}{e^{t^*}-1} = 2-\kappa.  \nonumber
    \end{align}
    The function $\frac{t}{e^t-1}$ is strictly decreasing from 1 to 0. Since $2-\kappa \in (0,1)$, a unique solution $t^*$ exists. Numerical evaluation of the function for $1<\kappa<2$ confirms that the global maximum is around $0.84$, which is less than 1. So, we can conclude that $F^{\mbox{exp}}_{\omega_1}(\phi)>F^{\mbox{par}}_{\omega_1}(\phi)$ since $R(t)<1$.    
\end{IEEEproof}

\section{Derivation of PDF of the Satellite Elevation Angle}\label{appen:elevation_pdf}

In this appendix, we derive the exact PDF of the satellite elevation angle $\theta$ observed by a ground user. We consider a geometric model involving the Earth's center ($O$), the user on the Earth's surface ($A$), and a satellite ($S$) on a sphere with radius $R_{orb}=R_E+h_{sat}$.

Let $\Gamma$ be the random variable representing the geocentric angle $\angle AOS$. Since satellites are uniformly distributed on the sphere, the probability of a satellite residing in a differential spherical zone at angle $\gamma$ is proportional to the surface area, which scales with $\sin \gamma$. Thus, the PDF of $\Gamma$ is given by:
\begin{align}
    f_{\Gamma}(\gamma) = \frac{\sin \gamma}{1 - k}, \quad 0 \le \gamma \le \arccos(k),\nonumber
\end{align}
where $k = R_E / R_{orb}$.

To find the PDF of the elevation angle $\Theta$, we use:
\begin{align}\label{eq:transform}
    f_{\Theta}(\theta) = f_{\Gamma}(\gamma) \left| \frac{\mathrm{d}\gamma}{\mathrm{d}\theta} \right|.
\end{align}
Let $d(\theta)$ denote the slant range (distance between $A$ and $S$). The relationship between $\theta$, $\gamma$, and $d(\theta)$ is governed by the sine and cosine laws in $\triangle OAS$:
\begin{align}
    \frac{R_{orb}}{\cos \theta} &= \frac{d(\theta)}{\sin \gamma} \implies \sin \gamma = \frac{d(\theta) \cos \theta}{R_{orb}}, \label{eq:sine_law} \\
    d(\theta)^2 &= R_E^2 + R_{orb}^2 - 2R_E R_{orb} \cos \gamma. \label{eq:cosine_law_gamma}
\end{align}
Differentiating \eqref{eq:cosine_law_gamma} with respect to $\theta$ involves the derivative $d'(\theta) \triangleq \frac{\mathrm{d}}{\mathrm{d}\theta}d(\theta)$. This yields:
{\begin{align}\label{eq:last}
2d(\theta) d'(\theta) = 2 R_E R_{orb} \sin \gamma \frac{\mathrm{d}\gamma}{\mathrm{d}\theta} .
\end{align}}
Substituting $\sin \gamma$ from \eqref{eq:sine_law} into \eqref{eq:last}, the Jacobian becomes
\begin{align}
    \frac{\mathrm{d}\gamma}{\mathrm{d}\theta} = \frac{1}{R_E \cos \theta} d'(\theta).\nonumber
\end{align}
Next, from the fact that $d(\theta)^2 + R_E^2 + 2R_E d(\theta) \sin \theta = R_{orb}^2$, we find $d'(\theta)$ by implicit differentiation:
\begin{align}
    d'(\theta) = - \frac{R_E d(\theta) \cos \theta}{d(\theta) + R_E \sin \theta}.\nonumber
\end{align}
Substituting these into \eqref{eq:transform}, we obtain
\begin{align}
    f_{\Theta}(\theta) &= \frac{1}{1-k} {\left( \frac{d(\theta) \cos \theta}{R_{orb}} \right)} {\left| \frac{1}{R_E \cos \theta} \left( - \frac{R_E d(\theta) \cos \theta}{d(\theta) + R_E \sin \theta} \right) \right|}\nonumber\\
    &= \frac{d(\theta)^2}{R_{orb}(1-k)(d(\theta) + R_E \sin \theta)} \cos \theta.\nonumber
\end{align}
Finally, substituting the explicit form of $d(\theta) = \sqrt{R_{orb}^2 - R_E^2 \cos^2 \theta} - R_E \sin \theta$, the PDF of $\theta$ becomes
\begin{align}
    f_{\Theta}(\theta) = \frac{\cos \theta}{(1-k)\sqrt{1 - k^2 \cos^2 \theta}} \left( \sqrt{1 - k^2 \cos^2 \theta} - k \sin \theta \right)^2.\nonumber
\end{align}

\section*{Acknowledgements}
The work of J. Lee was supported by the National Research Foundation of Korea(NRF) grant funded by the Korea government(MSIT) (RS-2026-25473096).

The work of F. Baccelli was supported by the Horizon Europe INSTINCT project (grant SNS 101139161), the France 2030 PEPR réseaux du Futur project (grant ANR-22-PEFT-0010), and by 5G NTN mmWave (BPIFrance). This joint work was also supported by a South Korea–France Hubert Curien grant.

\bibliographystyle{ieeetran}
\bibliography{referenceBibs}

\end{document}